# WEIGHTED UNIFORM CONSISTENCY OF KERNEL DENSITY ESTIMATORS

By Evarist Giné[1], Vladimir Koltchinskii[2] and Joel Zinn[3]

*University of Connecticut, University of New Mexico and Texas A&M University*

Let $f_n$ denote a kernel density estimator of a continuous density $f$ in $d$ dimensions, bounded and positive. Let $\Psi(t)$ be a positive continuous function such that $\|\Psi f^\beta\|_\infty < \infty$ for some $0 < \beta < 1/2$. Under natural smoothness conditions, necessary and sufficient conditions for the sequence $\sqrt{\frac{nh_n^d}{2|\log h_n^d|}}\|\Psi(t)(f_n(t) - Ef_n(t))\|_\infty$ to be stochastically bounded and to converge a.s. to a constant are obtained. Also, the case of larger values of $\beta$ is studied where a similar sequence with a different norming converges a.s. either to 0 or to $+\infty$, depending on convergence or divergence of a certain integral involving the tail probabilities of $\Psi(X)$. The results apply as well to some discontinuous not strictly positive densities.

**1. Introduction.** Over forty years ago, Parzen (1962) studied basic properties of kernel density estimators following their introduction by Rosenblatt (1956). Since then the kernel density estimator has become a classical object looked at by both statisticians and probabilists. For statisticians, it has been a canonical example of nonparametric curve estimators, which brought many important ideas from approximation theory and harmonic analysis into nonparametric statistics. Probabilists used the study of this estimator to test the strength of the methods from weak and strong convergence, empirical processes and probability in Banach spaces. In this paper, we consider a couple of problems about asymptotic behavior of kernel density estimators uniformly over all of $\mathbf{R}^d$ that do not seem to have been considered before,

Received January 2002; revised December 2003.

[1] Supported in part by NSF Grant DMS-00-70382.

[2] Supported in part by NSA Grant MDA904-99-1-0031.

[3] Supported in part by NSA Grant MDA904-01-1-0027.

particularly in the 1980s, when the basic results on uniform a.s. convergence were obtained.

The kernel density estimator $f_n$ of $f$ corresponding to a sample of size $n$, a kernel $K$ and a bandwidth $h > 0$ is

$$(1.1) \qquad f_n(t) = \frac{1}{nh^d} \sum_{i=1}^{n} K\left(\frac{X_i - t}{h}\right),$$

where $X_i$ are i.i.d. with density $f$. To ensure its consistency, $h$ is chosen to be a function $h_n$ of $n$ such that $h_n \to 0$ and $nh_n \to \infty$ as $n \to \infty$. This is a biased estimator, but we will not deal with the bias; we will only be interested in the sup norm of the deviation of $f_n$ from its mean.

Our starting point is the following well-known result due to Stute (1984):

$$(1.2) \qquad \lim_{n\to\infty} \sqrt{\frac{nh_n^d}{2|\log h_n^d|}} \left\| \frac{f_n - Ef_n}{\sqrt{f}} \right\|_J = \|K\|_2 \qquad \text{a.s.,}$$

where $J$ is a compact parallellepiped with sides parallel to the axes, $\| \cdot \|_J$ means "sup over $J$," $f$ is a uniformly continuous density which is bounded away from 0 on $J$, and $K$ is continuous and satisfies some additional assumptions [see, e.g., condition (K)]. Much later it was shown that

$$(1.3) \qquad \lim_{n\to\infty} \sqrt{\frac{nh_n^d}{2|\log h_n^d|}} \|f_n - Ef_n\|_\infty = \|K\|_2 \|f\|_\infty^{1/2} \qquad \text{a.s.,}$$

where $K$ satisfies condition (K) and $f$ is uniformly continuous [Giné and Guillou (2002) for any $d$, and Deheuvels (2000) for $d = 1$; a weaker result of this type was obtained much earlier by Silverman (1978)]. In both results the bandsequence $\{h_n\}$ satisfies Stute's (1982) conditions. In fact, these results can be slightly extended as follows: if $\Psi$ is uniformly continuous and bounded on $\bar{J}$, where $J$ is either a bounded parallellepiped of $\mathbf{R}^d$ with sides parallel to the axes, or $J = \mathbf{R}^d$, then

$$(1.4) \qquad \lim_{n\to\infty} \sqrt{\frac{nh_n^d}{2|\log h_n^d|}} \|\Psi(t)(f_n(t) - Ef_n(t))\|_J = \|K\|_2 \|\Psi f^{1/2}\|_J \qquad \text{a.s.,}$$

a result formulated in Deheuvels (2000) for $d = 1$ and which follows for any $d$ from Einmahl and Mason (2000) and Giné and Guillou (2002) (with simple modifications in their proofs). Note that (1.4) contains (1.2) and (1.3).

The first question on which we wish to shed some light is whether one can interpolate between the two results (1.2) and (1.3) by replacing $J$ by $\mathbf{R}^d$ and $f^{-1/2}$ by $f^{-\beta}$ for some $0 < \beta \le 1/2$ in (1.2). A more general formulation of the same problem is whether unbounded functions $\Psi$ are allowed in (1.4) when $J = \mathbf{R}^d$.

Notice that, in case $f > 0$ over all of $\mathbf{R}^d$ and $\liminf_{|x| \to \infty} f(x) = 0$, (1.4) implies that *only powers of $\beta$ not exceeding $1/2$ can lead to finite a.s. limits for the sequence*

$$(1.5) \qquad \left\{ \sqrt{\frac{nh_n^d}{2|\log h_n^d|}} \left\| \frac{f_n - Ef_n}{f^\beta} \right\|_\infty \right\}_{n=1}^\infty.$$

This is the case of classical norming, and in this case we find necessary and sufficient conditions (on the density $f$ and on the bandsequence $h_n$) for (1.5) to be stochastically bounded (Theorem 2.1); in fact, Theorem 2.1 gives necessary and sufficient conditions for

$$(1.5') \qquad \left\{ \sqrt{\frac{nh_n^d}{2|\log h_n^d|}} \| \Psi(t)(f_n(t) - Ef_n(t)) \|_\infty \right\}_{n=1}^\infty$$

to be stochastically bounded, assuming $\| \Psi f^\beta \|_\infty < \infty$ for some $\beta \in (0, 1/2)$. This result further clarifies the role of the sequence of maximum terms $\max_{1 \leq i \leq n} \Psi(X_i) / \sqrt{nh_n^d |\log h_n|}$ in the asymptotic behavior of (1.5′) in probability or in law. We also obtain a necessary and sufficient condition for (1.5′) to converge a.s. to the constant $\| K \|_2 \| \Psi f^{1/2} \|_\infty$ and show that if this condition is violated, then the sequence (1.5′) is a.s. unbounded (Theorem 2.6).

A second question is that of determining the right norming constants in the sequences (1.5) or (1.5′) for larger values of $\beta$ in order to obtain convergence. In this case, we also give necessary and sufficient conditions for stochastic boundedness (Theorem 3.1) and for a.s. convergence of the sequences (Theorem 3.4). The almost sure limit is shown to be either 0 or $+\infty$, depending on convergence or divergence of a certain integral describing the tail behavior of $\Psi(X)$. The situation in this case is somewhat similar to what is well known about weighted empirical processes; see Einmahl and Mason (1985a, b, 1988).

We consider a slightly more general situation where $f$ need not be strictly positive, however, we still require that, if $B_f = \{f > 0\}$, then $f$ be bounded away from zero on $B_f \cap \{|t| \leq a\}$ for all $a > 0$. Even this case requires unusual but somewhat natural smoothness conditions on $f$. More general situations seem to require a strengthening of the smoothness conditions, and we refrain here from considering them (see, however, Example 2.12).

*Assumptions and notation.* We introduce here some notation and conditions that are used throughout the paper.

For $x = (x_1, \ldots, x_d) \in \mathbf{R}^d$, we set $|x| := \max_{1 \leq i \leq d} |x_i|$. We assume that the kernel $K$ satisfies the following condition:

(K) $K \geq 0$, $K \not\equiv 0$, is a bounded measurable function with support contained in $[-1/2, 1/2]^d$ which belongs to the linear span (the set of finite linear combinations) of functions $k \geq 0$ satisfying the following property: the subgraph of $k$, $\{(s, u) : k(s) \geq u\}$, can be represented as a finite number of Boolean operations among sets of the form

$$\{(s, u) : p(s, u) \geq \varphi(u)\},$$

where $p$ is a polynomial on $\mathbf{R}^d \times \mathbf{R}$ and $\varphi$ is an arbitrary real function.

Conditions of a similar type were used, for example, in Koltchinskii and Sakhanenko (2000).

In particular, the above property is satisfied if the subgraph of $k$ is a semialgebraic set in $\mathbf{R}^d \times \mathbf{R}$ [see Dudley (1999), page 165]. If $K(x) = \phi(p(x))$, $p$ being a polynomial and $\phi$ a real function of bounded variation, then $K$ satisfies (K) [see Nolan and Pollard (1987)].

Condition (K) is mainly imposed because if $K$ satisfies it, then the class of functions

$$\mathcal{F} = \left\{ K\left(\frac{\cdot - t}{h}\right) : t \in \mathbf{R}^d, h > 0 \right\}$$

has covering numbers

$$N(\mathcal{F}, L_2(P), \|K\|_{L_2(P)}\varepsilon) \leq \left(\frac{A}{\varepsilon}\right)^v, \qquad 0 < \varepsilon < 1,$$

for some $A$ and $v$ finite and positive and for all probability measures $P$. Indeed, for a fixed polynomial $p$, the family of sets

$$\{\{(s, u) : p((s - t)/h, u) \geq \varphi(u)\} : t \in \mathbf{R}^d, h > 0\}$$

is contained in the family of positivity sets of a finite-dimensional space of functions, and then the entropy bound follows by Theorems 4.2.1 and 4.2.4 in Dudley (1999). The entropy bound will be crucial in the proofs below. Since the map $(x, t, h) \mapsto (x - t)/h$ is jointly measurable and $K$ is measurable, the class $\mathcal{F}$ is image admissible Suslin [Dudley (1999), page 186], and this implies that the measurability of the empirical process indexed by $\mathcal{F}$ [or even by $\{\Psi(t)K((\cdot - t)/h)\}$ with $\Psi$ continuous] is as good as if the class were countable; that is, we can ignore measurability of the sup of the empirical process over $g \in \mathcal{F}$ [cf. Dudley (1999), Corollary 5.3.5 and Theorem 5.3.6, or Pollard (1984), pages 195–197]. We set $\kappa := \|K\|_\infty$ (which is strictly positive).

The following assumptions on the density $f$ will be used repeatedly:

(D.a) $f$ is a bounded density on $\mathbf{R}^d$ continuous on its positivity set $B_f := \{t \in \mathbf{R}^d : f(t) > 0\}$, which is assumed to be open, and $\lim_{a \to \infty} \sup_{|t| > a} f(t) = 0$.

(D.b) For all $\delta > 0$, there exist $c \in (0, \infty)$ and $h_0 > 0$ such that, for all $|y| \le h_0$ and all $x \in B_f$, $x + y \in B_f$,

$$\frac{1}{c} f^\delta(x) \le \frac{f(x+y)}{f(x)} \le c f^{-\delta}(x).$$

(D.c) For all $r > 0$,

$$\lim_{h \to 0} \sup_{\substack{x, y: \, f(x) \ge h^r \\ x+y \in B_f, |y| \le h}} \left| \frac{f(x+y)}{f(x)} - 1 \right| = 0.$$

In particular, if $\log f$ is uniformly continuous on $\mathbf{R}^d$, then conditions (D.a)–(D.c) are satisfied (this is true, e.g., for the symmetric exponential density or for uniformly continuous nonvanishing densities with power tails). The above conditions are satisfied as well by normal and double exponential densities even though their logarithms are not uniformly continuous. Note also that (D.b) implies $\inf_{x \in B_f, |x| < a} f(x) > 0$ for all $a < \infty$ such that $B_f \cap \{|x| < a\} \ne \varnothing$, in particular, a continuous density with bounded support does not satisfy (D.b). Similarly, a density that has an isolated zero where it is continuous does not satisfy condition (D.b) either. In fact, Example 2.12 shows that, for such a density, the stochastic boundedness of the sequence (1.5) depends on the local behavior of the density at its zero points, and not only on the tails of the random variable $f^{-\beta}(X)$, as is the case under condition (D.b) (see Theorem 2.1). On the other hand, the exponential density does satisfy (D.a)–(D.c).

Conditions (D.b) and (D.c) on $f$ are not found in Stute's result (1982, 1984) or in Einmahl and Mason (2000) because they consider bounded intervals with $f$ bounded away from zero on them, and they are not found either in Giné and Guillou (2002) since there is no division by a power of $f$ in their result. These conditions seem natural for the results that will follow and we will indicate below that conditions of this type are indeed needed; see Example 2.11.

We assume that the weight function $\Psi$ satisfies the following conditions that resemble the above conditions on the density:

(W.a) $\Psi : B_f \mapsto \mathbf{R}_+$ is a positive continuous function on $B_f$.

(W.b) For all $\delta > 0$, there exist $c \in (0, \infty)$ and $h_0 > 0$ such that, for all $|y| \le h_0$ and all $x \in B_f$, $x + y \in B_f$,

$$\frac{1}{c} \Psi^{-\delta}(x) \le \frac{\Psi(x+y)}{\Psi(x)} \le c \Psi^\delta(x).$$

(W.c) For all $r > 0$,

$$\lim_{h \to 0} \sup_{\substack{x, y: \, \Psi(x) \le h^{-r} \\ x+y \in B_f, |y| \le h}} \left| \frac{\Psi(x+y)}{\Psi(x)} - 1 \right| = 0.$$

In particular, by (W.b), $\Psi$ is bounded on bounded subsets of $B_f$, but $\Psi$ may be unbounded if $B_f$ is unbounded.

We also need the following conditions that establish a relationship between $f$ and $\Psi$:

(WD.a)$_\beta$  $\|f^\beta\|_{\Psi,\infty} := \sup_{t \in B_f} |\Psi(t) f^\beta(t)| < \infty$, where $\beta$ is a positive number.

(WD.b)  For all $r > 0$,

$$\lim_{h \to 0} \sup_{\substack{x,y:\, \Psi(x) \le h^{-r} \\ x, x+\tau \in B_f, |y| \le h}} \left| \frac{f(x+y)}{f(x)} - 1 \right| = 0.$$

Note that (WD.a)$_\beta$ and (WD.b) imply (D.c): if $\|\Psi f^\beta\|_{B_f} \le c$ and $f(t) \ge h^r$, then $\Psi(t) \le ch^{-r\beta}$.

Also, if $\Psi \equiv f^{-\beta}$ (which is our main example), then (WD.a)$_\beta$ is satisfied and the set of conditions (D.a)–(D.c) is equivalent to (W.a)–(W.c) and (WD.b).

Regarding the window sizes, the assumptions are:

(H1)  $h_t$, $t \ge 1$, is monotonically decreasing to 0 and $t h_t^d$ is a strictly increasing function diverging to infinity as $t \to \infty$, and

(H2)  $h_t^d$ is regularly varying at infinity with exponent $-\alpha$ for some $\alpha \in (0,1)$; in particular, there exist $0 < \eta_0 \le \eta_1 < 1$ such that

$$\limsup_{t \to \infty} t^{\eta_0} h_t^d = 0 \quad \text{and} \quad \liminf_{t \to \infty} t^{\eta_1} h_t^d = \infty.$$

Condition (H2) is quite restrictive compared to the bandsequence assumptions in Stute ([1982](#)): besides the extra regularity, we do not allow $h_t$ to get too close to the extremes $1/t$ or $1/\log t$, and in particular, $|\log h_t|$ is comparable to $\log t$, $t > 1$. If we set

$$\lambda_t = \sqrt{t h_t^d |\log h_t|},$$

then, under (H1) and (H2), the function $\lambda_t$ is strictly increasing and is regularly varying with exponent larger than 0. This property of $\lambda_t$ is used throughout Section 2.

Our results rely on the by now classical theorem of Stute ([1984](#)) about the a.s. behavior of the uniform deviation of the kernel density estimator over compact intervals, suitably modified. The version of his theorem we need is a reformulation along the lines of Deheuvels ([2000](#)) of Proposition 3.1 in Giné and Guillou ([2002](#)), which in turn is adapted from Einmahl and Mason ([2000](#)).

PROPOSITION 1.1. *Let $f$ be a density on $\mathbf{R}^d$, continuous on an open set containing $D_a := \{t : |t| \le a, f(t) \ge a^{-1}\}$, for some $0 < a < \infty$. Let $\Psi$ be a strictly positive function, continuous on an open set containing $D_a$. Then*

$$(1.6) \qquad \lim_{n \to \infty} \sqrt{\frac{nh_n^d}{2|\log h_n^d|}} \|\Psi(t)(f_n - Ef_n)(t)\|_{D_a} = \|K\|_2 \|\Psi f^{1/2}\|_{D_a} \qquad a.s.$$

We omit the proof as it coincides with the proof of the abovementioned proposition, except for obvious changes.

Proposition 1.1 applies to $f$ satisfying (D.a) and (D.b) and $\Psi$ satisfying (W.a) and (W.b).

Without further mentioning, *all the results we state in this paper beyond this point assume conditions* (K), (H1), (H2), (D.a)–(D.c), (W.a)–(W.c), (WD.b) *and* (WD.a)$_\beta$ *for some $\beta$. The number $\beta$ is to be specified at each instance.* We will refer to these assumptions as the "usual hypotheses."

Finally, we introduce the following notation, which will be used throughout: for any function $g$ defined on $B_f$, we set

$$(1.7) \qquad \|g\|_{\Psi,\infty} := \sup_{t \in B_f} |g(t)\Psi(t)|.$$

## 2. The classical norming case.

The following theorem describes the stochastic boundedness behavior of the sequence (2.1). It shows in particular that no interpolation between (1.2) and (1.3) works for all strictly positive, bounded, continuous densities, and that when it works, it does not work for all the range of possible bandsequences. In what follows, $X$ is a random variable with density $f$.

THEOREM 2.1. *Assume the usual hypotheses, with condition* (WD.a)$_\beta$ *holding for some $\beta \in (0, 1/2)$, and, moreover, that either $B_f = \mathbf{R}^d$ or $K(0) = \|K\|_\infty$. Then the sequence*

$$(2.1) \qquad \left\{ \sqrt{\frac{nh_n^d}{|\log h_n|}} \|f_n - Ef_n\|_{\Psi,\infty} \right\}_{n=1}^\infty$$

*is stochastically bounded if and only if*

$$(2.2) \qquad \limsup_{t \to \infty} t \Pr\{\Psi(X) > (th_t^d|\log h_t|)^{1/2}\} < \infty.$$

*Moreover, under condition* (2.2), *the sequence*

$$(2.3) \qquad \begin{aligned} &\left\{ \sqrt{\frac{nh_n^d}{2|\log h_n^d|}} \|f_n - Ef_n\|_{\Psi,\infty} \right. \\ &\left. - \left( \max_{1 \le i \le n} \frac{\|K\|_\infty \Psi(X_i)}{\sqrt{2nh_n^d|\log h_n^d|}} \right) \vee (\|K\|_2\|f^{1/2}\|_{\Psi,\infty}) \right\}_{n=1}^\infty \end{aligned}$$

*converges to zero in probability.*

PROOF.  We will use the notation $\lambda_t = (th_t^d |\log h_t|)^{1/2}$. As mentioned above, conditions (H1) and (H2) imply that $\lambda_t$ is regularly varying with strictly positive exponent. Note that, by regular variation, condition (2.2) is equivalent to

$$(2.4) \qquad \limsup_{t \to \infty} t \Pr\{\Psi(X) > c\lambda_t\} < \infty$$

for any $0 < c < \infty$. By Montgomery-Smith's (1993) maximal inequality [see, e.g., de la Peña and Giné (1999), page 6] the stochastic boundedness of the sequence (2.1) implies that of the sequence

$$\left\{ \max_{1 \le i \le n} \frac{\|K((X_i - t)/h_n) - EK((X - t)/h_n)\|_{\Psi,\infty}}{\lambda_n} \right\}.$$

Then, since, for all $t$, letting $u = (u_1, \ldots, u_d)$,

$$EK((X - t)/h_n) = h_n^d \int_{-1/2}^{1/2} \cdots \int_{-1/2}^{1/2} K(u) f(h_n u + t) \, du_1 \cdots du_d$$

$$\le h_n^d \|K\|_1 \|f\|_\infty \to 0,$$

taking $t = X_1 - \tau h_n, \ldots, X_n - \tau h_n$ for $\tau \in \mathbf{R}^d$ satisfying $K(\tau) > 0$, we obtain that the sequence

$$\left\{ \max_{1 \le i \le n} \frac{\Psi(X_i - \tau h_n) I(X_i - \tau h_n \in B_f)}{\lambda_n} \right\}$$

is stochastically bounded. We take $\tau = 0$ if $K(0) \ne 0$. Now, if $K(0) \ne 0$, then $X_i - \tau h_n = X_i \in B_f$ a.s., and if $B_f = \mathbf{R}^d$, then obviously $X_i - \tau h_n \in B_f$, so that, in either case, $\Pr\{X - \tau h_n \in B_f\} = 1$. Thus, the sequence

$$\left\{ \max_{1 \le i \le n} \frac{\Psi(X_i - \tau h_n)}{\lambda_n} \right\}$$

is stochastically bounded. In particular, if $\tau = 0$,

$$\left\{ \max_{1 \le i \le n} \frac{\Psi(X_i)}{\lambda_n} \right\}$$

is stochastically bounded, proving condition (2.4) in this case (as, if $\xi_i$ are independent, $\Pr\{\max |\xi_i| > c\} \ge \sum \Pr\{|\xi_i| > c\}/(1 + \sum \Pr\{|\xi_i| > c\})$. If $\tau \ne 0$ but $B_f = \mathbf{R}^d$, given $\varepsilon > 0$, let $M$ be such that

$$\sup_n \Pr\left\{ \max_{1 \le i \le n} \frac{\Psi(X_i - \tau h_n)}{\lambda_n} > M \right\} < \varepsilon.$$

If $\Psi(X_i - \tau h_n) \le M\lambda_n$, then by regular variation there exists $r > 0$ such that $\Psi(X_i - \tau h_n) \le h_n^{-r}$ (at least for all $n$ large enough), and we can apply condition (W.c) to conclude that there exists $c > 1$ such that, for all $n$ large enough (independent of $X_i$), $\Psi(X_i) \leq c\Psi(X_i - \tau h_n)$. Then, for these values of $n$, we obtain

$$\Pr\left\{\max_{1 \leq i \leq n} \frac{\Psi(X_i)}{\lambda_n} > cM\right\} < \varepsilon.$$

Therefore, in this case, the sequence $\{\max_{1 \leq i \leq n} \Psi(X_i)/\lambda_n\}$ is also stochastically bounded, proving (2.4).

For the converse we note first that Proposition 1.1 takes care of the sup over $D_a$ for any $a > 0$.

Next, we observe that the centering in (2.1) can be ignored for a certain range of $t$'s. Let $\varepsilon_n \to 0$ and $0 < \delta < 1 - \beta$. Choose $r > 0$ such that

$$\frac{h_n^{r(1-\beta-\delta)}}{(nh_n^d)^{-1}\lambda_n} \to 0.$$

Then there exist $c < +\infty$ and $n_0 < \infty$ such that, by (D.b) and (D.c), for any $t \in B_f$ and $n \geq n_0$,

$$\frac{n\Psi(t)EK((X-t)/h_n)}{\lambda_n}$$

$$\leq \frac{c\kappa n h_n^d}{\lambda_n}\Psi(t) \sup_{\substack{|u| \leq 1/2 \\ t+h_n u \in B_f}} f(t + h_n u)$$

$$\leq \frac{c\kappa n h_n^d}{\lambda_n}\Psi(t)f(t)I(f(t) > h_n^r) + \frac{c\kappa n h_n^d}{\lambda_n}\Psi(t)f^{1-\delta}(t)I(f(t) \leq h_n^r).$$

Since, by condition (WD.a)$_\beta$, $\Psi(t) \leq cf^{-\beta}(t)$, the last summand tends to 0 uniformly in $t$ (as $\beta < 1 - \delta$). The sup of the first summand over all $t$ such that $\Psi(t)f(t) \leq \varepsilon_n^{1-\beta}((nh_n^d)^{-1}|\log h_n|)^{1/2}$ tends to 0 as well. Thus, we can ignore the centering $Ef_n(t)$ for all $t \in B_f$ such that

$$(2.5) \qquad \Psi(t)f(t) \leq \varepsilon_n^{1-\beta}\left(\frac{|\log h_n|}{nh_n^d}\right)^{1/2}$$

for any sequence $\varepsilon_n \to 0$. We take $\varepsilon_n = 1/\log n$.

In the rest of the proof, we consider the sup of $|\Psi(t)(f_n - Ef_n)(t)|$ over several regions.

First, we consider the regions

$$(2.6) \qquad A_n := \{t \in B_f : \Psi(t) > c_n^\beta \lambda_n\}$$

for the sequence $c_n = (\lambda_{n\log n}/\lambda_n)^{1/\beta}$ which tends to infinity because $\lambda_t$ is regularly varying with positive exponent. Actually, if $\eta > 0$ is the exponent of regular variation of $\lambda_t$, the representation formula for regularly varying functions [e.g., Feller ([1971](#)), page 282] gives that, for every $0 < \varepsilon < \eta$ and $c > 1$, there exists $n_0 < \infty$ such that

$$(2.7) \qquad \frac{1}{c} (\log n)^{\eta - \varepsilon} \leq \frac{\lambda_{n \log n}}{\lambda_n} \leq c (\log n)^{\eta + \varepsilon}$$

for all $n \geq n_0$. Then, since $\beta < 1/2 < 1 - \beta$, for a suitable $\delta > 0$ and all $n$ large enough, we have

$$c_n \lambda_n^{1/\beta} \geq n^{\eta/\beta - \delta} \geq n^{\eta/(1-\beta) + \delta}$$
$$\geq \frac{1}{\varepsilon_n} \left( \frac{\lambda_n}{|\log h_n|} \right)^{1/(1-\beta)} = \frac{1}{\varepsilon_n} \left( \frac{n h_n^d}{|\log h_n|} \right)^{1/(2(1-\beta))}.$$

This yields, for all $t \in A_n$,

$$\Psi(t) \geq \frac{1}{\varepsilon_n^\beta} \left( \frac{n h_n^d}{|\log h_n|} \right)^{\beta/(2(1-\beta))}$$

and, using condition (WD.a)$_\beta$ (which without loss of generality can be written as $\|f^\beta\|_{\Psi,\infty} \leq 1$), we get

$$f(t) \leq \varepsilon_n \left( \frac{|\log h_n|}{n h_n^d} \right)^{1/2(1-\beta)}.$$

This implies (2.5) for all $t \in A_n$, since $\Psi(t) f(t) \leq \|f^\beta\|_{\Psi,\infty} f^{1-\beta}(t) \leq f^{1-\beta}(t)$ [again due to (WD.a)$_\beta$]. Therefore,

$$\sup_{t \in A_n} \frac{n \Psi(t) E K((X-t)/h_n)}{\lambda_n} \to 0,$$

showing that we can ignore the centering $E f_n$ on the region $A_n$. For any point $a = (a_1, \ldots, a_d) \in \mathbf{R}^d$ and positive number $\rho$, we set

$$J(a; \rho) := [a_1 - \rho/2, a_1 + \rho/2] \times \cdots \times [a_d - \rho/2, a_d + \rho/2] \cap B_f.$$

Then, discarding the centering,

$$\sup_{t \in A_n} \frac{\Psi(t) \sum_{i=1}^n K((X_i - t)/h_n)}{\lambda_n} \leq \sup_{t \in A_n} \frac{\kappa}{\lambda_n} \Psi(t) \sum_{i=1}^n I(X_i \in J(t; h_n)).$$

Now we divide $A_n$ into two parts:

$$A_{n,1} := \{ t \in B_f : \Psi(t) > h_n^{-r} \} \quad \text{and} \quad A_{n,2} := \{ t \in B_f : c_n^\beta \lambda_n < \Psi(t) \leq h_n^{-r} \},$$

where $r$ is such that $h_n^{-r(1-\delta)} \geq c_n^\beta \lambda_n$ for some $\delta > 0$ and all $n$. It follows from condition (W.b) that $t \in A_{n,1}$ and $s \in J(t; h_n)$ imply that there are $c$ and $n_0$ such that

$$\frac{\Psi(s)}{\Psi(t)} \geq c \Psi^{-\delta}(t)$$

for all $n \geq n_0$, so that $\Psi(s) \geq c\Psi^{1-\delta}(t) > ch_n^{-r(1-\delta)}$. Hence, for the same values of $n$ and some $C < \infty$, we have

$$\Pr\left\{\sup_{t \in A_{n,1}} \frac{\Psi(t)\sum_{i=1}^n K((X_i - t)/h_n)}{\lambda_n} > \varepsilon\right\} \leq n \Pr\{\Psi(X) \geq h_n^{-r(1-\delta)}\}$$

$$\leq \frac{C}{\log n} \to 0.$$

It follows from condition (W.c) that $t \in A_{n,2}$ and $s \in J(t; h_n) \cap B_f$ imply that there are $c$ and $n_0$ such that $\Psi(s) \geq c_n^\beta \lambda_n / c$, for all $n \geq n_0$. Hence, for these values of $n$ we have

$$\Pr\left\{\sup_{t \in A_{n,2}} \frac{\Psi(t)\sum_{i=1}^n K((X_i - t)/h_n)}{\lambda_n} > \varepsilon\right\} \leq n \Pr\left\{\Psi(X) \geq \frac{c_n^\beta \lambda_n}{c}\right\} \leq \frac{C}{\log n} \to 0$$

for some $C < \infty$ by (2.4). The last two limits imply that

$$(2.8) \qquad \lim_{n \to \infty} \sup_{t \in A_n} \frac{\Psi(t)\sum_{i=1}^n K((X_i - t)/h_n)}{\lambda_n} = 0 \qquad \text{in pr.}$$

Now we consider the regions

$$(2.9) \qquad B_n := \left\{t \in B_f : f(t)\Psi(t) \leq \varepsilon_n^{1-\beta}\left(\frac{|\log h_n|}{nh_n^d}\right)^{1/2},\right.$$

$$\left. \Psi(t) \leq c_n^\beta (nh_n^d |\log h_n|)^{1/2}\right\}$$

and notice that in these regions we can also ignore the centering [by (2.5)]. Our goal is to show that

$$\left\{\sup_{t \in B_n} \frac{\Psi(t)\sum_{i=1}^n K((X_i - t)/h_n)}{\lambda_n}\right\}_{n=1}^\infty$$

is stochastically bounded under condition (2.2) and that, moreover, if either $B_f = \mathbf{R}^d$ or $K(0) = \kappa$, then also

$$(2.10) \qquad \sup_{t \in B_n} \frac{\Psi(t)\sum_{i=1}^n K((X_i - t)/h_n)}{\lambda_n} = \kappa \max_{1 \leq i \leq n} \frac{\Psi(X_i)}{\lambda_n} + o_p(1).$$

As above,

$$\sup_{t \in B_n} \frac{\Psi(t)\sum_{i=1}^n K((X_i - t)/h_n)}{\lambda_n} \leq \sup_{t \in B_n} \frac{\kappa\Psi(t)}{\lambda_n}\sum_{i=1}^n I(X_i \in J(t; h_n)),$$

and we set

$$Z_n := \sup_{t \in B_n} \frac{\Psi(t)}{\lambda_n}\sum_{i=1}^n I(X_i \in J(t; h_n)).$$

For $j = 1, \ldots, n$, set

$$B_{n,j} := B_n \cap J(X_j; h_n).$$

If $t \notin \bigcup_{j=1}^n B_{n,j}$, then $Z_n = 0$. Hence, we have

$$Z_n = \max_{1 \le j \le n} \sup_{t \in B_{n,j}} \frac{\Psi(t)}{\lambda_n} \sum_{i=1}^n I(X_i \in J(t; h_n)).$$

By conditions (W.c) and (WD.b), $t \in B_{n,j}$ implies that

$$\Psi(t) \le c\Psi(X_j)$$

and also

$$\Psi(X_j) \le cc_n^\beta \lambda_n, \qquad f(X_j)\Psi(X_j) \le c\varepsilon_n^{1-\beta}\left(\frac{|\log h_n|}{nh_n^d}\right)^{1/2}$$

for any $c > 1$, provided that $n$ is large enough.

Set

$$I_j = I_{n,j} := I\left(\Psi(X_j) \le cc_n^\beta \lambda_n, f(X_j)\Psi(X_j) \le c\varepsilon_n^{1-\beta}\left(\frac{|\log h_n|}{nh_n^d}\right)^{1/2}\right).$$

Then

$$
\begin{aligned}
(2.11) \quad Z_n &\le \max_{1 \le j \le n} \frac{c\Psi(X_j)I_j \sum_{i=1}^n I(|X_i - X_j| \le h_n)}{\lambda_n} \\
&\le \max_{1 \le j \le n} \frac{c\Psi(X_j)}{\lambda_n} + \max_{1 \le j \le n} \frac{c\Psi(X_j)I_j \sum_{1 \le i \le n, i \ne j} I(|X_i - X_j| \le h_n)}{\lambda_n}.
\end{aligned}
$$

By condition (2.2), the first term in the above bound is the general term of a stochastically bounded sequence. We now show that the second term tends to zero in probability. To handle this term, let $P_j$ denote conditional expectation given $X_j$ and set

$$p_j := P_j\{|X - X_j| \le h_n\}.$$

It follows from condition (D.a) that

$$2^d c^{-1} h_n^d f(X_j) \le p_j \le 2^d c h_n^d f(X_j)$$

(provided that $I_j = 1$). A standard bound on binomial probabilities [e.g., Giné and Zinn (1984), page 958] shows that

$$
\begin{aligned}
P_j&\left\{\frac{I_j\Psi(X_j)\sum_{1 \le i \le n, i \ne j} I(|X_i - X_j| \le h_n)}{\lambda_n} \ge \varepsilon\right\} \\
&\le \left(\frac{(n-1)ep_j\Psi(X_j)}{\lambda_n\varepsilon}\right)^{(\varepsilon\lambda_n/\Psi(X_j))\vee 1}.
\end{aligned}
$$

Using the bound on $p_j$, this probability can be further bounded by

$$\left(\frac{2^d e c n h_n^d f(X_j) \Psi(X_j)}{\lambda_n \varepsilon}\right)^{(\varepsilon \lambda_n / \Psi(X_j)) \vee 1}.$$

We can and do assume that $I_j = 1$ (otherwise the conditional probability in question is 0). Then

$$f(X_j) \Psi(X_j) \leq c \varepsilon_n^{1-\beta} \left(\frac{|\log h_n|}{n h_n^d}\right)^{1/2},$$

and we have

$$\frac{2^d e c n h_n^d f(X_j) \Psi(X_j)}{\lambda_n \varepsilon} \leq \frac{C_1 n h_n^d \varepsilon_n^{1-\beta} (|\log h_n| / n h_n^d)^{1/2}}{(n h_n^d |\log h_n|)^{1/2} \varepsilon} = \frac{C_1 \varepsilon_n^{1-\beta}}{\varepsilon}$$

for some $C_1 < \infty$ (and all $n$ large enough). Note also that

$$\frac{2^d e c n h_n^d f(X_j) \Psi(X_j)}{\lambda_n \varepsilon} = \frac{C}{\varepsilon} \left(\frac{n h_n^d}{|\log h_n|}\right)^{1/2} f(X_j) \Psi(X_j),$$

where $C$ is a finite positive cosntant. For large $n$, $\frac{C_1 \varepsilon_n^{1-\beta}}{\varepsilon} \leq e^{-1/\varepsilon}$, which yields

$$P_j \left\{ \frac{I_j \Psi(X_j) \sum_{1 \leq i \leq n, i \neq j} I(|X_i - X_j| \leq h_n)}{\lambda_n} \geq \varepsilon \right\}$$

$$\leq \left(\exp\left\{-\frac{\lambda_n}{\Psi(X_j)}\right\}\right) \wedge \left(\frac{C}{\varepsilon} \left(\frac{n h_n^d}{|\log h_n|}\right)^{1/2} f(X_j) \Psi(X_j)\right).$$

Let

$$I_j^1 := I\left(\Psi(X_j) \leq \frac{\lambda_n}{(3 \log n)}, f(X_j) \Psi(X_j) \leq c \varepsilon_n^{1-\beta} \left(\frac{|\log h_n|}{n h_n^d}\right)^{1/2}\right)$$

and let $I_j^2 := I_j - I_j^1$. Then we have

$$\Pr\left\{\max_{1 \leq j \leq n} \frac{I_j \Psi(X_j) \sum_{1 \leq i \leq n, i \neq j} I(|X_i - X_j| \leq h_n)}{\lambda_n} \geq \varepsilon\right\}$$

$$\leq \sum_{j=1}^n E I_j P_j \left\{\frac{I_j \Psi(X_j) \sum_{1 \leq i \leq n, i \neq j} I(|X_i - X_j| \leq h_n)}{\lambda_n} \geq \varepsilon\right\}$$

$$\leq \sum_{j=1}^n E I_j^1 \exp\left\{-\frac{\lambda_n}{\Psi(X_j)}\right\} + \sum_{j=1}^n E I_j^2 \frac{C}{\varepsilon} \left(\frac{n h_n^d}{|\log h_n|}\right)^{1/2} f(X_j) \Psi(X_j)$$

$$=: (I) + (II).$$

Then, using the definition of $I_j^1$ and $I_j^2$, we get

$$(I) \leq n \exp\left\{-\frac{3\log n}{\lambda_n}\lambda_n\right\} = n^{-2}$$

and

$$(II) \leq n\frac{C}{\varepsilon}\left(\frac{nh_n^d}{|\log h_n|}\right)^{1/2}\frac{(3\log n)^{(1-\beta)/\beta}}{(nh_n^d|\log h_n|)^{(1-\beta)/(2\beta)}}\Pr\left\{\Psi(X) \geq \left(\frac{\lambda_n}{3\log n}\right)\right\}.$$

Now, since $\lambda_t$ is regularly varying with a strictly positive exponent, the representation theorem for regularly varying functions gives that $\lambda_n/(3\log n) \geq c\lambda_{n/(\log n)^\gamma}$ for some $\gamma > 0$, $c > 0$ and all $n$ large enough [see (2.7)]. Hence, by (2.4), there exists $C > 0$ such that, for these values of $n$,

$$(II) \leq \frac{C}{\varepsilon}\frac{(\log n)^{(1-\beta)/\beta+\gamma}}{|\log h_n|^{1/2+(1-\beta)/(2\beta)}}\frac{1}{(nh_n^d)^{(1-\beta)/(2\beta)-1/2}}.$$

By (H1) and (H2), this is at most of the order of logarithmic factors times $n^{(1-\eta_1)[1/2-(1-\beta)/(2\beta)]}$, a negative power of $n$ because $0 < \beta < 1/2$. Thus, $(II)$ also tends to zero. Since both $(I)$ and $(II)$ tend to 0, we have

$$\Pr\left\{\max_{1\leq j\leq n}\frac{I_j\Psi(X_j)\sum_{1\leq i\leq n, i\neq j}I(|X_i-X_j|\leq h_n)}{\lambda_n} \geq \varepsilon\right\} \to 0 \qquad \text{as } n\to\infty.$$

This implies [see bound (2.11)] that

$$(2.12)\quad \sup_{t\in B_n}\frac{\Psi(t)\sum_{i=1}^n K((X_i-t)/h_n)}{\lambda_n} \leq \kappa Z_n \leq c\kappa\max_{1\leq i\leq n}\frac{\Psi(X_i)}{\lambda_n} + o_p(1),$$

for any $c > 1$. The stochastic boundedness of

$$\left\{\sup_{t\in B_n}\frac{\Psi(t)\sum_{i=1}^n K((X_i-t)/h_n)}{\lambda_n}\right\}$$

follows immediately from this inequality and condition (2.2).

To bound the supremum from below, choose $\tau$ such that $K(\tau) > \kappa - \delta$ (for a small $\delta$) with the understanding that if $K(0) = \kappa$, then we choose $\tau = 0$, so that either $\tau = 0$ or $B_f = \mathbf{R}^d$. Then

$$\sup_{t\in B_n}\frac{\Psi(t)\sum_{i=1}^n K((X_i-t)/h_n)}{\lambda_n} \geq (\kappa-\delta)\max_{1\leq i\leq n}\frac{\Psi(X_i-\tau h_n)I_{B_n}(X_i-\tau h_n)}{\lambda_n}.$$

Hence, in view of this and the two-sided bound for $\kappa Z_n$ immediately above, to establish (2.10), it is enough to show that

$$(2.13)\quad \max_{1\leq i\leq n}\frac{\Psi(X_i-\tau h_n)I_{B_n}(X_i-\tau h_n)}{\lambda_n} = \max_{1\leq i\leq n}\frac{\Psi(X_i)}{\lambda_n} + o_p(1).$$

Since condition (W.c) implies that, for any $c > 1$ and for large enough $n$,

$$c^{-1} < \frac{\Psi(X_i - \tau h_n)}{\Psi(X_i)} < c$$

(assuming that $X_i - \tau h_n \in B_n$), taking $c$ arbitrarily close to 1 reduces the proof of (2.13) to showing that

$$\max_{1 \leq i \leq n} \frac{\Psi(X_i) I_{B_n}(X_i - \tau h_n)}{\lambda_n} = \max_{1 \leq i \leq n} \frac{\Psi(X_i)}{\lambda_n} + o_p(1),$$

or, put in another way, (2.13) will be proved if we show that

$$(2.13') \qquad \max_{1 \leq i \leq n} \frac{\Psi(X_i) I_{B_n^c}(X_i - \tau h_n)}{\lambda_n} \to 0 \qquad \text{in pr.}$$

$B_n^c$ naturally decomposes into the union of three regions and we look separately at each of them. If $B_f = \mathbf{R}^d$, then $I_{B_f^c}(X_i - \tau h_n) = 0$, and if $\tau = 0$, then this indicator is 0 a.s., so that, in either case,

$$\max_{1 \leq i \leq n} \frac{\Psi(X_i) I_{B_f^c}(X_i - \tau h_n)}{\lambda_n} \to 0 \qquad \text{a.s.}$$

Next, we consider

$$\Pr\left\{ \max_{1 \leq i \leq n} \frac{\Psi(X_i) I(\Psi(X_i - \tau h_n) \geq c_n^\beta \lambda_n)}{\lambda_n} > \varepsilon \right\}$$

$$\leq n \Pr\{\Psi(X - \tau h_n) \geq c_n^\beta \lambda_n\}$$

$$\leq n \Pr\{c_n^\beta \lambda_n \leq \Psi(X - \tau h_n) \leq h_n^{-r}\} + n \Pr\{\Psi(X - \tau h_n) \geq h_n^{-r}\}.$$

Using condition (W.c), we get (for any $c > 1$)

$$n \Pr\{c_n^\beta \lambda_n \leq \Psi(X - \tau h_n) \leq h_n^{-r}\} \leq n \Pr\{\Psi(X) \geq c^{-1} c_n^\beta \lambda_n\} \to 0.$$

Similarly, using condition (W.b) (recall that $X - \tau h_n \in B_\Psi$ with probability 1),

$$n \Pr\{\Psi(X - \tau h_n) \geq h_n^{-r}\} \leq n \Pr\{\Psi(X) \geq c^{-1} h_n^{-r/(1+\delta)}\}$$

for some $c > 0$ and $\delta > 0$. Assuming that $r$ is large enough (so that $h_n^{-r/(1+\delta)} \geq c_n^\beta \lambda_n$), we then conclude that

$$n \Pr\{\Psi(X - \tau h_n) \geq h_n^{-r}\} \to 0$$

and, hence,

$$\Pr\left\{ \max_{1 \leq i \leq n} \frac{\Psi(X_i) I(\Psi(X_i - \tau h_n) \geq c_n^\beta \lambda_n)}{\lambda_n} > \varepsilon \right\} \to 0.$$

Before considering the last piece of $B_n^c$, we note that, since $f^\beta(t)\Psi(t) \le 1$ for all $t$, if moreover $f(u)\Psi(u) > L$, then $f^{1-\beta}(u) > L$ and consequently $\Psi(u) \le f^{-\beta}(u) < L^{-\beta/(1-\beta)}$, an observation that we will use several times below. This observation and condition (W.c) give

$$\max_{1 \le i \le n} \frac{\Psi(X_i)I(f(X_i - \tau h_n)\Psi(X_i - \tau h_n) > \varepsilon_n^{1-\beta}(|\log h_n|/(nh_n^d))^{1/2})}{\lambda_n}$$

$$\le \max_{1 \le i \le n} \frac{\Psi(X_i)I(\Psi(X_i - \tau h_n) < c\varepsilon_n^{-\beta}(nh_n^d/|\log h_n|)^{\beta/(2(1-\beta))})}{\lambda_n}$$

$$\le c\varepsilon_n^{-\beta}\left(\frac{nh_n^d}{|\log h_n|}\right)^{\beta/(2(1-\beta))}\frac{1}{(nh_n^d|\log h_n|)^{1/2}}$$

$$= \frac{c}{\varepsilon_n^\beta|\log h_n|^{1/2+\beta/(2(1-\beta))}(nh_n^d)^{1/2-\beta/(2(1-\beta))}}.$$

Now, since $\beta < 1/2 < 1 - \beta$ and $nh_n^d \ge n^{1-\eta_1}$ [by (H2)], whereas $\varepsilon_n = 1/\log n$ and $|\log h_n|$ is comparable to $\log n$, it follows that the above bound is dominated by a negative power of $n$ so that, in particular, it tends to zero. This and the previous two limits conclude the proof of $(2.13')$ and hence of (2.10).

Finally, we consider the sup over the remaining set of $t$'s. For $a$ large, fixed, just as above, set

$$(2.14) \quad C_n = C_{n,a} := \{t \in D_a^c \cap B_f : f(t)\Psi(t) \ge \varepsilon_n^{1-\beta}(|\log h_n|/nh_n^d)^{1/2}\},$$

where $\varepsilon_n$ is as defined in the previous paragraph. In this range the centering cannot be ignored. We will apply an estimate for the expected supremum of the empirical process over bounded Vapnik–Červonenkis type classes of functions [Giné and Guillou ([2001](#)), inequality (2.1) and Talagrand ([1994](#)), for classes of sets; see also Einmahl and Mason (2000) for a similar inequality]: if a class of functions $\mathcal{F}$ is measurable (in particular, if it is image admissible Suslin) and satisfies

$$(2.15) \qquad N(\mathcal{F}, L_2(Q), \varepsilon\|F\|_\infty) \le \left(\frac{A}{\varepsilon}\right)^v, \qquad 0 < \varepsilon < 1,$$

for some $v \ge 1$, $A \ge 3\sqrt{e}$ finite and all finite probability measures $Q$, where $F$ is a measurable envelope for the class $\mathcal{F}$, then

$$(2.16) \qquad E\|n(P_n - P)\|_\mathcal{F} \le C\left(\sqrt{v}\sqrt{n}\sigma\sqrt{\log\frac{AU}{\sigma}} + vU\log\frac{AU}{\sigma}\right),$$

where $\sigma$ and $U$ are any numbers satisfying $0 < \sigma < U$ and

$$(2.17) \qquad \sigma^2 \ge \sup_{g \in \mathcal{F}}\mathrm{Var}_P(g), \qquad U \ge \sup_{g \in \mathcal{F}}\|g\|_\infty,$$

and $C$ is a universal constant. [In Giné and Guillou (2001), condition (2.15) has $\|F\|_{L_2(Q)}$ instead of $\|F\|_\infty$, but it can be easily checked that their proof works as well under condition (2.15).] As mentioned immediately below the statement of condition (K), there exist $A$ and $v$ finite such that

$$N(\{K((\cdot - t)/h_n) : t \in \mathbf{R}\}, L_2(Q), \varepsilon) \leq \left(\frac{A\kappa}{\varepsilon}\right)^v, \qquad 0 < \varepsilon < 1,$$

for all $h_n > 0$ and all probability measures $Q$ on $\mathbf{R}$. Now, the class of functions

$$\mathcal{F}_n := \{\Psi(t) K((\cdot - t)/h_n) : t \in C_n\}$$

is contained in

$$\mathcal{G}_n := \{u K((\cdot - t)/h_n) : t \in \mathbf{R}, 0 < u \leq U_n\},$$

where

$$(2.18) \qquad U_n := \frac{\kappa}{\varepsilon_n^\beta}\left(\frac{n h_n^d}{|\log h_n|}\right)^{\beta/(2(1-\beta))}$$

[recall that, as observed above, under condition (WD.a)$_\beta$, $f\Psi \geq \alpha$ implies $\Psi \leq \alpha^{-\beta/(1-\beta)}$]. Therefore, since the $L_2(Q)$ distance between $u K((\cdot - t)/h_n)$ and $v K((\cdot - s)/h_n)$ is dominated by $\kappa |u - v| + U_n\|K((\cdot - t)/h_n) - K((\cdot - s)/h_n)\|_{L_2(Q)}$, it follows by taking optimal coverings of $[0, U_n]$ with respect to the Euclidean distance, and of $\mathcal{F}_n$ with respect to the $L_2(Q)$ distance, that the entropy bound

$$(2.19) \qquad N(\mathcal{F}_n, L_2(Q), \varepsilon U_n) \leq \left(\frac{2A\kappa}{\varepsilon}\right)^{v+1}, \qquad 0 < \varepsilon < 1,$$

holds for all probability measures $Q$ and all $n$ large enough. The class $\mathcal{F}_n$ is also image admissible Suslin since the map $(x, t) \mapsto \Psi(t) K((x - t)/h_n)$ is measurable. So, inequality (2.16) applies to it. We can take $U = U_n$ as defined in (2.18). Next we estimate $\sigma_n^2$. It follows from a previous observation and from regular variation that, on $C_n$, we have both $f \geq h_n^r$ and $\Psi \leq h_n^r$ for some $r$ and all $n$ large enough. Then, (D.c) and (W.c) give that there exist $c, C, n_0 < \infty$ independent of $a$ such that, for all $n \geq n_0$ and all $t \in C_n = C_{n,a}$,

$$\Psi^2(t) E K^2((X - t)/h_n) \leq c E(K^2((X - t)/h_n)\Psi^2(X))$$

$$= c h_n^d \int_{\substack{|u| \leq 1/2 \\ t + h_n u \in B_f}} K^2(u)\Psi^2(t + h_n u) f(t + h_n u)\, du$$

$$\leq c h_n^d \|K\|_2^2 \|f\Psi^2\|_{D_a^c \cap B_f}$$

$$\leq C h_n^d (\|f\Psi^2\|_{D_a^c \cap B_f} \vee n^{-1}).$$

So, we can take

$$\sigma_n^2 = C h_n^d (\|f\Psi^2\|_{D_a^c \cap B_f} \vee n^{-1}).$$

The constant $A = A_n$ must be taken to be $(2A\kappa) \vee (3\sqrt{e})$, where $A$ is the constant in (2.15) for the class consisting of translations and dilations of $K$. In particular, since by (H2) $|\log h_n|$ is comparable to $\log n$, we have

$$\log \frac{A_n U_n}{\sigma_n} \le c |\log h_n|$$

for some constant $c < \infty$ independent of $n$. So, inequality (2.16) applied to $\mathcal{F}_n$ gives

$$E \sup_{t \in C_n} \left| \frac{\Psi(t) \sum_{i=1}^n (K((X_i-t)/h_n) - EK((X-t)/h_n))}{\lambda_n} \right|$$

$$\le \frac{C}{\lambda_n} \left[ \lambda_n (\|f^{1/2}\Psi\|_{D_a^c \cap B_f} \vee n^{-1/2}) + \frac{1}{\varepsilon_n^\beta} \left( \frac{n h_n^d}{|\log h_n|} \right)^{\beta/(2(1-\beta))} |\log h_n| \right]$$

for a constant $C$ independent of $n$, for all sufficiently large $n$. We should note that the numerical constants in the above inequalities are not only independent of $n$, but they are independent of $a$ as well. Since $\beta < 1/2$ and therefore $\beta/(1-\beta) < 1$, and since, by (D.a) and (WD.a)$_\beta$,

$$\|f^{1/2}\Psi\|_{D_a^c \cap B_f} \le \|f^{1/2-\beta}\|_{D_a^c \cap B_f} = \|f\|_{D_a^c \cap B_f}^{1/2-\beta} \to 0 \qquad \text{as } a \to \infty,$$

we obtain

$$(2.20) \quad \begin{aligned} &\lim_{a \to \infty} \limsup_{n \to \infty} E \sup_{t \in C_{n,a}} \left| \frac{\Psi(t) \sum_{i=1}^n (K((X_i-t)/h_n) - EK((X-t)/h_n))}{\lambda_n} \right| \\ &\le \lim_{a \to \infty} C \|f^{1/2}\Psi\|_{D_a^c \cap B_f} = 0. \end{aligned}$$

Now, the theorem follows from (1.6), (2.8), (2.10) and (2.20). $\square$

Now we make two comments on the assumptions.

REMARK 2.2.   The assumption "$B_f = \mathbf{R}^d$ or $K(0) = \|K\|_\infty$" has been imposed because in general we may not have $X - \tau h_n \in B_f$ with small enough probability, as $n \Pr\{X - \tau h_n \notin B_f\}$ could well be of the order of $n h_n \to \infty$. Now, this condition has been used in full only in the proof of (2.3). Proving that tightness of the sequence (2.1) implies condition (2.2) has only required $B_f = \mathbf{R}^d$ or $K(0) > 0$, whereas proving that condition (2.2) implies tightness of the sequence (2.1) does not require any hypothesis of this type.

The above proof justifies, a posteriori, having taken $\beta < 1/2$:

COROLLARY 2.3. *Assume* (K), (H1), (H2), (D.a)–(D.c) *and* $B_f = \mathbf{R}^d$. *Then the sequence*

$$\left\{ \sqrt{\frac{nh_n^d}{2|\log h_n^d|}} \left\| \frac{f_n - Ef_n}{\sqrt{f}} \right\|_\infty \right\}_{n=1}^\infty$$

[*which coincides with* (2.1) *for* $\Psi = f^{-1/2}$] *is not stochastically bounded.*

PROOF. By the first part of the previous proof, if (2.1) with $\Psi = f^{-1/2}$ is tight, then there is $C > 0$ such that

$$n \Pr\left\{ \frac{1}{f(X)} > \lambda_n^2 \right\} \leq C.$$

Since $f$ takes all the values between 0 and $\|f\|_\infty$, for $n$ large enough there is $x_n$ in $\mathbf{R}^d$ such that $f(x_n) = 1/(2\lambda_n^2)$. Then, by condition (D.c), there is a subset $D_n$ containing $x_n$ and of Lebesgue measure larger than $\lambda_n^{-2/r}$, where $1/f(x) \geq \lambda_n^2$ and $f(x) \geq 1/(4\lambda_n^2)$, and therefore, if we take $r \geq \eta_1/d$ with $\eta_1$ as in condition (H2),

$$n \Pr\left\{ \frac{1}{f(X)} > \lambda_n^2 \right\} \geq n \Pr\{X \in D_n\} \geq \frac{n}{4\lambda_n^{2(1+d/r)}} \to \infty,$$

contradiction. □

Theorem 2.1 has the following obvious corollary regarding convergence in distribution:

COROLLARY 2.4. *Under the assumptions in Theorem* 2.1, *the sequence* (2.1) *converges in distribution if and only if the sequence of maxima,*

$$\left\{ \max_{1 \leq i \leq n} \frac{\Psi(X_i)}{\lambda_n} \right\},$$

*converges in distribution. Then, if $Z$ is a random variable with distribution the limit of this last sequence, we have*

$$\sqrt{\frac{nh_n^d}{2\log h_n^{-d}}} \|f_n - Ef_n\|_{\Psi,\infty} \xrightarrow{d} (\|K\|_\infty Z) \vee (\|K\|_2 \|f^{1/2}\|_{\Psi,\infty}).$$

Next we consider the a.s. counterpart to Theorem 2.1. The following proposition will help. It is perhaps relevant to recall first a well-known fact, whose proof we omit as it is similar to a classical result of Feller [e.g., Lemma 3.2.4, Corollary 3.2.3 and Theorem 3.2.5 in Stout (1974)].

LEMMA 2.5.  *Let $V_i$ be i.i.d. real random variables and let $\{c(n)\}$ be a nondecreasing sequence, regularly varying with strictly positive exponent. Then,*

$$either \qquad \limsup_{n \to \infty} \max_{1 \le i \le n} \frac{|V_i|}{c(n)} = \infty \qquad a.s. \quad or \quad \lim_{n \to \infty} \max_{1 \le i \le n} \frac{|V_i|}{c(n)} = 0 \qquad a.s.$$

*And this happens according to whether*

$$\sum_n \Pr\{|V_n| > Cc(n)\} = \infty \quad or \quad \sum_n \Pr\{|V_n| > Cc(n)\} < \infty$$

*for some (or, equivalently, all) $C > 0$.*

PROPOSITION 2.6.  *Assume that conditions* (D.a), (W.a)–(W.c) *and* (WD.a)$_\beta$, (WD.b) *hold for some $\beta > 0$ and that, moreover, either $B_f = \mathbf{R}^d$ or $K(0) > 0$. Let $c(n) \nearrow \infty$ be a regularly varying function of $n$. Assume*

$$(2.21) \quad \limsup_n \left\| \frac{\sum_{i=1}^n (K((X_i - t)/h_n) - EK((X - t)/h_n))}{c(n)} \right\|_{\Psi, \infty} < \infty \qquad a.s.$$

*Then,*

$$(2.22) \qquad \sum_n \Pr\left\{ \frac{\Psi(X)}{c(n)} > C \right\} < \infty$$

*for all $0 < C < \infty$ or, what turns out to be the same by Lemma 2.5,*

$$(2.22') \qquad \lim_{n \to \infty} \max_{1 \le i \le n} \frac{\Psi(X_i)}{c(n)} = 0 \qquad a.s.$$

PROOF.  The proof is standard, but we give it here for completeness. First we note that if (2.21) holds and $c(n) \nearrow \infty$ is regularly varying, then $c(n)$ has necessarily positive exponent, which by Lemma 2.5 gives the equivalence between (2.22) and (2.22'). This follows because, by (2.21), there is $t$ with $\Psi(t) \ne 0$ and $f(t) \ne 0$ such that the sequence $\sum_{i=1}^n \Psi(t)(K(\frac{X_i - t}{h_n}) - EK(\frac{X - t}{h_n}))/c(n)$, $n \in \mathbf{N}$, is tight, which, by boundedness and finite support of $K$, implies that the sequence of its second moments is uniformly bounded, thus, that the sequence $nh_n^d/c^2(n)$ is bounded; hence, since by (H2) $nh_n^d$ is regularly varying with strictly positive exponent and $c(n)$ is regularly varying, it follows that the exponent of $c(n)$ is strictly positive as well.

Let $\{X_i'\}$ be an independent copy of $\{X_i\}$. We can symmetrize in (2.21) and still have the lim sup finite. By continuity of $\Psi$ on $B_f$, there is $n(\omega) < \infty$ a.s. such that, for all $n \ge n(\omega)$,

$$\left\| \frac{K((X(\omega) - t)/h_n) - K((X'(\omega) - t)/h_n)}{c(n)} \right\|_{\Psi, \infty} \le \frac{\kappa(\Psi(X) + \Psi(X') + 1)}{c(n)}.$$

This tends to zero and therefore the lim sup in (2.21) is a.s. constant by the zero–one law. Hence, we have

$$\Pr\left\{\sup_{n \geq k}\left\|\frac{\sum_{i=1}^{n}\Psi(t)(K((X_i-t)/h_n) - K((X_i'-t)/h_n))}{c(n)}\right\|_{B_f} > c\right\} \to 0$$
$$\text{as } k \to \infty$$

for some $c < \infty$. Set

$$H_n(X, X') := \frac{\Psi(t)(K((X-t)/h_n) - K((X'-t)/h_n))}{c(n)},$$

and, for $k \in \mathbf{N}$,

$$Z_{i,k} = (H_k(X_i, X_i'), H_{k+1}(X_i, X_i'), \dots, H_{k+r}(X_i, X_i'), \dots)$$

if $i \leq k$, and

$$Z_{i,k} = (0, \overset{r)}{\dots}, 0, H_{k+r}(X_{k+r}, X_{k+r}'), H_{k+r+1}(X_{k+r}, X_{k+r}') \dots)$$

for $i = k + r$, $r = 1, \dots$. Then, the above sup over $n \geq k$ is simply

$$\left\|\sum_{i=1}^{\infty} Z_{i,k}\right\|,$$

where $\|(x_1(t), \dots, x_n(t), \dots)\| = \sup_n \|x_n(t)\|_{B_f}$. The random vectors $Z_{i,k}$ are independent and symmetric, and we can apply Lévy's inequality to get that

$$\Pr\left\{\sup_{i \in \mathbf{N}} \|Z_{i,k}\| > 2c\right\} \to 0$$

as $k \to \infty$. By independence, this implies that

$$\sum_{i=1}^{\infty} \Pr\{\|Z_{i,k}\| > 2c\} \to 0$$

as $k \to \infty$. Let $\tau = 0$ if $K(0) > 0$ and otherwise let $|\tau| < 1$ be such that $K(\tau) > 0$. Then,

$$\|H_m(X, X')\|_{B_f} \geq \frac{\tilde{c}\Psi(X - h_m\tau)}{c(n)} I(|X - X'| > h_n)$$

for some $\tilde{c} > 0$, and we get that

$$\|Z_{i,k}\| = \sup_{m \geq i \vee k} H_m(X_i, X_i') \geq \frac{\tilde{c}\Psi(X_i)}{c(i \vee k)} I(|X_i - X_i'| > h_{i \vee k})$$

when $\tau = 0$ and

$$\|Z_{i,k}\| = \sup_{m \geq i \vee k} H_m(X_i, X_i') \geq \frac{\tilde{c}\inf_{|h| \leq |h_{i \vee k}|} \Psi(X_i - h)}{c(i \vee k)} I(|X_i - X_i'| > h_{i \vee k})$$

when $B_f = \mathbf{R}^d$. The case $\tau = 0$ is easier to handle, so we will complete the proof only for the second case. In this case, since $\Pr'\{|X - X'| > h_i\} \geq 1 - \|f\|_\infty h_i^d$, the previous inequality yields

$$\sum_{i=1}^\infty \Pr\{\|Z_{i,k}\| > 2c\} \geq \sum_{i \geq k} (1 - \|f\|_\infty h_i^d) \Pr\left\{ \frac{\tilde{c} \inf_{|h| \leq |h_i|} \Psi(X - h)}{c(i)} > 2c \right\}.$$

Then by (W.b) there are $0 < \delta < 1$ and $\hat{c} > 0$ such that

$$\sum_{n=1}^\infty \Pr\{\Psi(X) > \hat{c} c^{1/(1+\delta)}(n)\} < \infty.$$

But by regular variation, there exists $r > 0$ such that $h_n^{-r} > \hat{c} c^{1/(1+\delta)}(n)$, and therefore

$$\sum_{n=1}^\infty \Pr\{\Psi(X) > h_n^{-r}\} < \infty.$$

Now by (W.c), for $n$ large enough, there exists $C < \infty$ such that

$$\Pr\{\Psi(X) > Cc(n)\} \leq \Pr\left\{ \frac{\tilde{c} \inf_{|h| \leq |h_n|} \Psi(X - h)}{c(n)} > 2c \right\} + \Pr\{\Psi(X) > h_n^{-r}\}.$$

Therefore,

$$\sum_{n=1}^\infty \Pr\{\Psi(X) > Cc(n)\} < \infty. \qquad \square$$

We are now prepared to give an integral test for a.s. convergence of the sequence (2.1). Notice the difference with the tightness criterion, which is due to the fact that, by Lemma 2.5, we have

(2.23)  either  $\displaystyle\lim_{n\to\infty} \max_{1 \leq i \leq n} \frac{\Psi(X_i)}{\lambda_n} = 0$ a.s.  or  $\displaystyle\limsup_{n\to\infty} \max_{1 \leq i \leq n} \frac{\Psi(X_i)}{\lambda_n} = \infty$.

THEOREM 2.7. *Assume the usual hypotheses, with condition* (WD.a)$_\beta$ *holding for some* $\beta \in (0, 1/2)$, *and, moreover, that either* $B_f = \mathbf{R}^d$ *or* $K(0) = \|K\|_\infty$. *Set* $\lambda(t) = \sqrt{t h_t^d |\log h_t|}$, *as before. Then, either*

$$(2.24) \qquad \lim_{n\to\infty} \sqrt{\frac{n h_n^d}{2|\log h_n^d|}} \|f_n - Ef_n\|_{\Psi,\infty} = \|K\|_2 \|f^{1/2}\|_{\Psi,\infty} \qquad a.s.$$

*or*

$$(2.25) \qquad \limsup_{n\to\infty} \sqrt{\frac{n h_n^d}{2|\log h_n^d|}} \|f_n - Ef_n\|_{\Psi,\infty} = \infty \qquad a.s.,$$

*according to whether*

$$(2.26) \quad \int_1^\infty \Pr\{\Psi(X) > \lambda_t\}\, dt < \infty \quad or \quad \int_1^\infty \Pr\{\Psi(X) > \lambda_t\}\, dt = \infty.$$

PROOF. By Proposition 2.6, since $\lambda_n$ is regularly varying, if the integral in (2.26) is infinite, then (2.25) holds. So, we must prove that

$$(2.27) \qquad \int_1^\infty \Pr\{\Psi(X) > c\lambda_t\}\,dt < \infty$$

for all $c > 0$ implies (2.24). We proceed as in the proof of Theorem 2.1, with the addition of the usual blocking and replacing, in the estimation of the sup over $C_n$, the moment bound by an exponential inequality. By (2.4), we only have to consider the sup of our statistics over $A_n$, $B_n$ and $C_n$, the three sets defined as in the proof of Theorem 2.1, but with $c_n = 1$ (and $\varepsilon_n = 1/\log n$ as before), and we can ignore the centerings on $A_n$ and $B_n$. By monotonicity of $h_n$ and $\lambda_n$, we have

$$\max_{2^k \leq n \leq 2^{k+1}} \sup_{t \in A_n} \frac{\Psi(t)\sum_{i=1}^n K((X_i - t)/h_n)}{\lambda_n}$$

$$\leq \kappa \sup_{t \in A_{2^k}} \frac{\Psi(t)\sum_{i=1}^{2^{k+1}} I(X_i \in J(t; h_{2^k}))}{\lambda_{2^k}}.$$

Hence, we have, as before,

$$\Pr\left\{\max_{2^k \leq n \leq 2^{k+1}} \sup_{t \in A_n} \frac{\Psi(t)\sum_{i=1}^n K((X_i - t)/h_n)}{\lambda_n} > \varepsilon\right\}$$

$$\leq \varepsilon^{-1} 2^{k+1} \Pr\left\{\Psi(X) \geq \frac{\lambda_{2^k}^{1/\beta}}{c}\right\}$$

for all $k$ large enough and some $c > 0$. But, by (2.27), this is the general term of a convergent series, thus proving that

$$(2.28) \qquad \lim_{n \to \infty} \sup_{t \in A_n} \frac{\Psi(t)\sum_{i=1}^n K((X_i - t)/h_n)}{\lambda_n} = 0 \qquad \text{a.s.}$$

Regarding $B_n$ (with $c_n = 1$ and $\varepsilon_n \searrow 0$), we first note that, by regular variation,

$$\bigcup_{n=2^k}^{2^{k+1}} B_n \subseteq \tilde{B}_{2^k} := \left\{t : f(t)\Psi(t) \leq c'\varepsilon_{2^k}^{1-\beta}\left(\frac{|\log h_{2^k}|}{2^k h_{2^k}^d}\right)^{1/2},\right.$$

$$\left. \Psi(t) \leq c'c_{2^k}^\beta (2^k h_{2^k}^d |\log h_{2^k}|)^{1/2}\right\}$$

for some $c' > 1$. Then, as in (2.11),

$$\max_{2^k \leq n \leq 2^{k+1}} \sup_{t \in B_n} \frac{\Psi(t)\sum_{i=1}^n K((X_i - t)/h_n)}{\lambda_n}$$

$$\leq \sup_{t \in \tilde{B}_{2^k}} \frac{\Psi(t) \sum_{i=1}^{2^{k+1}} I(X_i \in J(t; h_{2^k}))}{\lambda_{2^k}}$$

$$\leq \max_{1 \leq j \leq 2^{k+1}} \frac{c\Psi(X_j)}{\lambda_{2^k}} + \max_{1 \leq j \leq 2^{k+1}} \frac{c\Psi(X_j)I_j \sum_{1 \leq i \leq n, i \neq j} I(|X_i - X_j| \leq h_{2^k})}{\lambda_{2^k}},$$

where $I_j$ is defined as before but with $n = 2^{k+1}$, and $c$ may be different from the constant in (2.11). Now, the maximum term tends to zero a.s. by (2.7) and Lemma 2.5, and the remainder term satisfies

$$\Pr\left\{ \max_{1 \leq j \leq 2^{k+1}} \frac{c\Psi(X_j)I_j \sum_{1 \leq i \leq n, i \neq j} I(|X_i - X_j| \leq h_{2^k})}{\lambda_{2^k}} > \varepsilon \right\} \leq \frac{C}{\varepsilon 2^{\alpha k}}$$

for some $\alpha > 0$ and all $k$ large enough, as in the proof of Theorem 2.1. Therefore,

$$(2.29) \qquad \lim_{n \to \infty} \sup_{t \in B_n} \frac{\Psi(t) \sum_{i=1}^{n} K((X_i - t)/h_n)}{\lambda_n} = 0 \qquad \text{a.s.}$$

In order to control the sup of our statistics over $C_n = C_{n,a}$ [as defined in (2.14)], we will use Talagrand's exponential inequality [Talagrand (1994, 1996)] in conjunction with the bound on the expected value of the sup of an empirical process given in (2.16). In a ready to use form for the problem at hand, it is as follows [Giné and Guillou (2001), equation (2.12)]: under assumption (2.15) above, and with the notation of (2.17) above, assuming further that

$$0 < \sigma < U/2 \quad \text{and} \quad \sqrt{n}\sigma \geq U\sqrt{\log \frac{U}{\sigma}},$$

there exist constants $C$ and $L$ such that, for all $s > C$,

$$(2.30) \Pr\left\{ \left\| \sum_{i=1}^{n} (f(\xi_i) - Ef(\xi_1)) \right\|_{\mathcal{F}} > s\sigma\sqrt{n}\sqrt{\log \frac{U}{\sigma}} \right\} \leq L\exp\left\{ -\frac{D(s)}{L}\log \frac{U}{\sigma} \right\},$$

where

$$D(s) := s\log(1 + s/4L) \to \infty \qquad \text{as } s \to \infty.$$

We apply this inequality to the class $\mathcal{F}_n$ defined on the last part of the proof of Theorem 2.1, with $U = U_n$ and $\sigma = \sigma_n$ as defined there, so that $\log \frac{U_n}{\sigma_n} \asymp \log n$. Since, for $a$ fixed and $n$ large enough, $\sigma_n \to 0$, $U_n \to \infty$ and $\sqrt{n}\sigma_n/(U_n\sqrt{\log \frac{U_n}{\sigma_n}}) \to \infty$, the above applies to give that there exists $C < \infty$ such that, for all $a > 0$ and for all $n$ large enough (depending on $a$),

$$(2.31) \quad \Pr\left\{ \sup_{t \in C_{n,a}} \left| \frac{\Psi(t) \sum_{i=1}^{n}(K((X_i - t)/h_n) - EK((X - t)/h_n))}{\lambda_n} \right| \right.$$
$$\left. > C \| f^{1/2} \Psi \|_{D_a^c \cap B_f} \right\}$$
$$\leq L \exp\{-2 \log n\}.$$

Hence,

$$(2.32) \quad \limsup_{n \to \infty} \sup_{t \in C_{n,a}} \left| \frac{\Psi(t) \sum_{i=1}^{n}(K((X_i - t)/h_n) - EK((X - t)/h_n))}{\lambda_n} \right|$$
$$\leq C \| f^{1/2} \Psi \|_{D_a^c \cap B_f} \quad \text{a.s.}$$

Combining (1.6), (2.28), (2.29) and (2.32), and letting $a \to \infty$, we obtain the limit (2.24). $\square$

We conclude this section with a few examples. We take $\Psi(t) = f^{-\beta}(t)$. Other choices of $\Psi$ are of course possible.

EXAMPLE 2.8.  Suppose $f : \mathbf{R} \mapsto (0, M]$ is continuous and

$$f(x) = c_1 e^{-c_2 |x|^r}$$

for all $|x|$ large enough, for some $r > 0$ and for some constants $c_1$ and $c_2$. Then, $f$ satisfies (D.a)–(D.c). Take

$$h_n = n^{-\alpha}, \qquad 0 < \alpha < 1.$$

For simplicity assume $c_1 = c_2 = 1$. It is easy to see that

$$\Pr\{|X| > u\} \asymp u^{1-r} e^{-u^r}.$$

Hence,

$$\Pr\left\{ \frac{1}{f(X)} > t^{(1-\alpha)/(2\beta)} (\log t)^{1/(2\beta)} \right\} \asymp \frac{1}{t^{(1-\alpha)/(2\beta)} (\log t)^{1/(2\beta)-(1-r)/r}}.$$

Then the above theorems imply the following. For $r \geq 1$, which includes the symmetric exponential and the normal densities, the conclusion is that the sequence (2.1) with $\Psi(t) = f^{-\beta}(t)$ is tight (stochastically bounded) if and only if

$$2\beta \leq 1 - \alpha$$

and that, if this is the case, then

$$(2.33) \quad \sqrt{\frac{n^{1-\alpha}}{2\alpha \log n}} \left\| \frac{f_n - Ef_n}{f^{\beta}} \right\|_{\infty} \to \|K\|_2 \|f\|_{\infty}^{1/2-\beta} \qquad \text{a.s.}$$

The same is true for exponential densities if we replace in (2.33) the sup over
**R** by the sup over $\mathbf{R}^+$. For $0 < r < 1$, if $2\beta < 1 - \alpha$, then the limit (2.33)
holds. If $2\beta = 1 - \alpha$, different behaviors arise; namely, if $(1 - r)/r > 1/(2\beta)$,
then the sequence (2.1) is not stochastically bounded; if $(1 - r)/r = 1/(2\beta)$,
the sequence converges in distribution to the limit in distribution of the random
variables

$$\left( \max_{1 \le i \le n} \frac{\|K\|_\infty}{\sqrt{2\alpha n^{2\beta} \log n} f^\beta(X_i)} \right) \vee (\|K\|_2 \|f\|_\infty^{1/2 - \beta}),$$

which is unbounded and can be easily computed (see the next example);
if $(1 - r)/r - 1/(2\beta) < 0$, we have convergence in probability in (2.33), but
convergence a.s. holds only if $(1 - r)/r - 1/(2\beta) < -1$.

EXAMPLE 2.9.    Suppose now the real density $f$ is strictly positive, continuous and

$$f(x) = \frac{c}{|x|^r}$$

for all $|x|$ large enough, for some $r > 1$ and for some constant $c$. These
densities also satisfy (D.a)–(D.c). Take $h_n = n^{-\alpha}$, $\alpha \in (0, 1)$ as above. Then,
(2.1) [again, with $\Psi(t) = f^{-\beta}(t)$] is tight if and only if

$$\beta \le \frac{r - 1}{r} \frac{1 - \alpha}{2},$$

and, if this is the case, then (2.33) holds true.

EXAMPLE 2.10.    Let now $f(x) = \frac{1}{2} e^{-|x|}$ be the symmetric exponential
density on **R**. Then,

$$\Pr\left\{ \max_{1 \le i \le n} \frac{1}{f(X_i)} > u \right\} = \begin{cases} 1 - \left( 1 - \dfrac{2}{u} \right)^n, & \text{if } u \ge 2, \\ 1, & \text{otherwise,} \end{cases}$$

so that

$$\max_{1 \le i \le n} \frac{1}{n^\beta f^\beta(X_i)} \xrightarrow{d} Z^\beta,$$

where $Z$ has distribution

$$\Pr\{Z \le t\} = \begin{cases} e^{-2/t}, & \text{if } t > 0, \\ 0, & \text{otherwise.} \end{cases}$$

Hence, if we take $\beta \in (0, 1/2)$ and

$$h_n = \frac{1}{n^{1 - 2\beta} \log n},$$

Theorem 2.1 gives that

$$\sqrt{\frac{nh_n}{2|\log h_n|}}\Big\|\frac{f_n - Ef_n}{f^\beta}\Big\|_\infty \xrightarrow{d} \max\Big(\frac{\|K\|_\infty}{\sqrt{2(1-2\beta)}}Z^\beta, \frac{\|K\|_2}{2^{1/2-\beta}}\Big).$$

The next two examples show that the above results are not true in general without conditions of the type of (D.b), (D.c) [and (W.b), (W.c)]. The first addresses smoothness and the second the existence of zeros of $f$ on the closure of $B_f$.

EXAMPLE 2.11. It is easy to see that the double exponential density still satisfies conditions (D.a)–(D.c) and, hence, Theorems 2.1 and 2.7, but the density

$$f(t) := ce^{-e^t}, \qquad t \geq 0,$$

does not. Specifically, condition (D.b) fails for this density and we show below that, for all $\beta \in (0,1)$ and for $h_n = n^{-\alpha}$,

$$(2.34) \qquad \sqrt{\frac{nh_n}{2|\log h_n|}}\sup_{t \geq 0}\Big|\frac{f_n(t) - Ef_n(t)}{f^\beta(t)}\Big| \to \infty \qquad \text{a.s.}$$

Indeed, if $K$ is continuous and strictly positive at the point $t = -1/4$, then

$$\frac{Ef_n(t)}{f^\beta(t)} = c\exp\{\beta e^t\}\frac{1}{h_n}EK\Big(\frac{X-t}{h_n}\Big)$$

$$= c\exp\{\beta e^t\}\int_{-1/2}^{1/2}K(u)f(h_n u + t)\,du$$

$$\geq c_1\exp\{\beta e^t - e^{t-4^{-1}n^{-\alpha}}\}$$

$$= c_1\exp\{e^t[\beta - e^{e^t(e^{-4^{-1}n^{-\alpha}}-1)}]\}.$$

Let $t_n := \log n$. Then, for large $n$,

$$\sqrt{\frac{nh_n}{2|\log h_n|}}\frac{Ef_n(t_n)}{f^\beta(t_n)} \geq c_1\exp\{e^n[\beta - e^{n(e^{-4^{-1}n^{-\alpha}}-1)}]\} \geq c_1\exp\Big\{\beta\frac{e^n}{2}\Big\}.$$

On the other hand,

$$\Pr\{f_n(t_n) \neq 0\} \leq \Pr\Big\{\max_{1 \leq i \leq n}X_i \geq t_n - \frac{n^{-\alpha}}{2}\Big\} \leq Cn\exp\{-e^{\sqrt{n}}\},$$

which implies that

$$\sqrt{\frac{nh_n}{2|\log h_n|}}\frac{f_n(t_n)}{f^\beta(t_n)} \to 0 \qquad \text{a.s.},$$

and therefore (2.34) holds.

EXAMPLE 2.12.   This example shows that if the density $f$ has a zero in **R**, then the asymptotic behavior of

$$\sqrt{\frac{nh_n}{|\log h_n|}}\left\|\frac{f_n - Ef_n}{f^\beta}\right\|_{B_f}$$

depends on the local behavior of $f$ at the zero point and is no longer controlled only by condition (2.2). Note that in this case condition (D.b) fails. For simplicity, assume that $h_n = n^{-\alpha}$ (with $\alpha < 1$) and $K = I_{[-1/2,1/2]}$. Let $f$ be a density continuous on a neighborhood of 0 and such that $f(0) = 0$ and, moreover, for some $s > 0$,

$$f(t) \asymp |t|^s \qquad \text{as } t \to 0.$$

In particular, we assume that $f$ is $s$ times continuously differentiable at 0 (for an even integer number $s$) and $f^{(j)}(0) = 0$ for $j < s$, $f^{(s)}(0) > 0$. It is easy to see that

$$(2.35) \qquad \Pr\{|X| \le t\} \asymp t^{s+1} \qquad \text{as } t \to 0.$$

We will show that if $s > \frac{1}{\alpha} - 1$, then, for all $C > 0$,

$$(2.36) \qquad \Pr\left\{\sqrt{\frac{nh_n}{2|\log h_n|}}\left\|\frac{f_n - Ef_n}{f^\beta}\right\|_{B_f} > C\right\} \to 1.$$

The proof is almost the same as in the previous example. Let $t_n \to 0$ be chosen in such a way that $f(t_n) = e^{-n}$. Note that $t_n = o(h_n)$. Then, using (2.35), we get

$$\sqrt{\frac{nh_n}{2|\log h_n|}}\frac{Ef_n(t_n)}{f^\beta(t_n)} \asymp \sqrt{\frac{n^{1-\alpha}}{\log n}}e^{\beta n}n^\alpha EK(n^\alpha(X - t_n))$$

$$\asymp \sqrt{\frac{n^{1+\alpha}}{\log n}}e^{\beta n}\Pr\left\{t_n - \frac{n^{-\alpha}}{2} \le X \le t_n + \frac{n^{-\alpha}}{2}\right\}$$

$$\asymp \sqrt{\frac{n^{1+\alpha}}{\log n}}e^{\beta n}n^{-(s+1)\alpha} \to \infty.$$

On the other hand, also using (2.35), if $s > \frac{1}{\alpha} - 1$, then

$$\Pr\{f_n(t_n) \ne 0\} \le \Pr\{\exists\, i, 1 \le i \le n : X_i \in (t_n - h_n/2, t_n + h_n/2)\}$$

$$\le n\Pr\{X \in (t_n - h_n/2, t_n + h_n/2)\}$$

$$\asymp nh_n^{s+1} = n^{1-\alpha(s+1)} \to 0.$$

This immediately implies (2.36). Now let $f(t) = c|t|^s$ for $|t| \le a$ and $f(t) = 0$ otherwise. Then it is easy to check that condition (2.2) holds if and only if $\beta \le \frac{1-\alpha}{2}(1 + \frac{1}{s})$. Thus, for large enough $s$, this condition does not imply the stochastic boundedness of (2.1).

**3. Large normings.** By Proposition 1.1, the central part of the process $\Psi(t)(f_n(t) - Ef_n(t))$, that is, its sup over $D_a$, for all $a > 0$, has an influence on the asymptotic size in probability of the sequence (2.1) and completely determines its a.s. limit. But if we normalize by a sequence larger than $\sqrt{nh_n^d|\log h_n|}$, this central part of the sup vanishes for all $a > 0$, and only the extremes of the range of $t$'s should have an influence on the limit. This is what we examine in this section. As in the previous section, we will only consider regularly varying window sizes and normings. As is to be expected, the only possible limit a.s. in this situation is zero, and the sum is asymptotically equivalent, in probability, to the maximum term. This is roughly the content of the following two theorems.

THEOREM 3.1. *Assume the usual hypotheses, with condition* (WD.a)$_\beta$ *holding for some $\beta \in (0,1]$, and, moreover, that either $B_f = \mathbf{R}^d$ or $K(0) = \|K\|_\infty$. Let $d_t$ be a strictly increasing regularly varying function such that $d_t/\lambda_t \to \infty$ and $d_t \geq Ct^\beta$ for some $C > 0$. Then, the sequence*

$$(3.1) \qquad \left\{ \left\| \frac{\sum_{i=1}^n (K((X_i - t)/h_n) - EK((X - t)/h_n))}{d_n} \right\|_{\Psi,\infty} \right\}$$

*is stochastically bounded if and only if*

$$(3.2) \qquad \limsup_{t \to \infty} t \Pr\{\Psi(X) > d_t\} < \infty.$$

*Moreover, if condition (3.2) holds, then*

$$(3.3) \qquad \begin{aligned} &\left\| \frac{\sum_{i=1}^n (K((X_i - t)/h_n) - EK((X - t)/h_n))}{d_n} \right\|_{\Psi,\infty} \\ &\quad - \max_{1 \leq i \leq n} \frac{\|K\|_\infty \Psi(X_i)}{d_n} \to 0 \qquad \text{in pr.} \end{aligned}$$

PROOF. The proof is similar to that of Theorem 2.1. First we consider $\beta < 1$. Necessity of condition (3.2) follows exactly in the same way. Here we indicate the few changes that should be made to the proof of Theorem 2.1 in order to prove that (3.2) implies (3.1) and (3.3). First, and this is by far the main difference with Theorem 2.1, the sup of

$$(3.4) \qquad \left| \frac{\Psi(t) \sum_{i=1}^n (K((X_i - t)/h_n) - EK((X - t)/h_n))}{d_n} \right|$$

over $D_a$, tends to zero a.s. for all $a < \infty$ by Proposition 1.1. Regarding the centering, consider the bound

$$\frac{n\Psi(t)EK((X - t)/h_n)}{d_n f^\beta(t)}$$

$$\leq \frac{c\kappa nh_n^d}{d_n}\Psi(t)f(t)I(f(t) > h_n^r) + \frac{c\kappa nh_n^d}{d_n}\Psi(t)f^{1-\delta}(t)I(f(t) \leq h_n^r),$$

where $t \in B_f$, $1 - \beta > \delta$ and $r$ is such that $n h_n^{d+r(1-\beta-\delta)}/d_n \to 0$, which is obtained as in the proof of Theorem 2.1. If the exponent of regular variation of $n h_n^d$ is strictly smaller than that of $d_n$, then, since $\Psi f^\eta$ is bounded for all $\eta \geq \beta$, the sup over $t \in B_f$ of this bound tends to zero and therefore we can simply ignore the centerings in (3.1) and (3.3). Otherwise, the second summand tends to zero uniformly in $t \in B_f$ and the first tends to zero uniformly on all $t \in B_f$ such that

$$f(t)\Psi(t) \leq \varepsilon_n^{1-\beta} \frac{d_n}{n h_n^d},$$

for any $\varepsilon_n \to 0$. So we can ignore the centerings for these values of $t$. As before, we take $\varepsilon_n = 1/\log n$.

Continuing in analogy with the proof of Theorem 2.1, we now define

$$A_n = \{t \in B_f : \Psi(t) > c_n^\beta d_n\}$$

with $c_n = (d_{n \log n}/d_n)^{1/\beta} \to \infty$, and we get, as in (2.8) but now using the properties of $d_n$, that

$$\lim_{n \to \infty} \sup_{t \in A_n} \frac{\sum_{i=1}^n \Psi(t) K((X_i - t)/h_n)}{d_n} = 0 \qquad \text{in pr.}$$

(for $0 < \beta < 1$).

Next we set

$$B_n := \left\{ t \in B_f : f(t)\Psi(t) \leq \varepsilon_n^{1-\beta} \frac{d_n}{n h_n^d}, \Psi(t) \leq c_n^\beta d_n \right\}$$

in analogy with (2.9). Then, proceeding as in the proof of (2.10) with the only formal change of replacing $\lambda_n$ by $d_n$ and $\sqrt{n h_n^d/|\log h_n|}$ by $n h_n^d/d_n$, we arrive at analogous conclusions, namely that the sequence

$$\sup_{t \in B_n} \frac{\Psi(t) \sum_{i=1}^n K((X_i - t)/h_n)}{d_n}$$

is stochastically bounded and that in fact it can be represented as

$$\max_{1 \leq i \leq n} \frac{\kappa \Psi(X_i)}{d_n} + o_p(1).$$

(This requires using the properties of $d_n$ and $h_n$ but, given that proof, the details are straightforward.)

Finally, we consider

$$C_n = B_f \setminus (A_n \cup B_n) = \{t \in B_f : f(t)\Psi(t) \geq \varepsilon_n^{1-\beta} d_n/(n h_n^d)\}.$$

Using as before that $\Psi f \geq L$ implies, by (WD.a)$_\beta$, that $\Psi \leq L^{-\beta/(1-\beta)}$, we can take

$$U_n = \kappa \varepsilon_n^{-\beta} \left( \frac{n h_n^d}{d_n} \right)^{\beta/(1-\beta)}.$$

We will consider two cases.

If the exponent of regular variation of $nh_n^d$ is strictly smaller than that of $d_n$, then $\varepsilon_n^{1-\beta} d_n/(nh_n^d) \to \infty$ and therefore, since, by (WD.a)$_\beta$, $\|f\Psi\|_\infty \le \|f\|_\infty^{1-\beta} < \infty$, $C_n$ is eventually the empty set.

Assume now that the exponent of regular variation of $d_n$ does not exceed that of $nh_n^d$. Then $\varepsilon_n^{1-\beta} d_n/(nh_n^d)$ is eventually dominated by $n^{-\delta}$ for any $\delta > 0$, so that we eventually have $f(t) \ge h_n^r$ and $\Psi(t) \le h_n^{-r}$ for some $r > 0$ and all $t \in C_n$. So, we can apply (D.c) and (W.c), which, together with (WD.a)$_\beta$, immediately imply that we can take $\sigma_n$ as follows:

$$\sigma_n^2 = \begin{cases} C\kappa h_n^d \|f\|_\infty^{1-2\beta}, & \text{if } \beta \le 1/2, \\ C\kappa h_n^d \varepsilon_n^{-(2\beta-1)} (nh_n^d/d_n)^{(2\beta-1)/(1-\beta)}, & \text{if } \beta > 1/2. \end{cases}$$

Since $U_n$ is either slowly varying or tends to infinity and $\sigma_n$ tends to zero as a negative power of $n$ for $\beta \le 1/2$, we get, in this case, that, eventually,

$$0 < \sigma_n < U_n/2 \quad \text{and} \quad \log \frac{U_n}{\sigma_n} \asymp \log n.$$

The same conclusion holds for $\beta > 1/2$ since $h_n^{d/2}$ decreases as a negative power of $n$ and the exponent of $nh_n^d/d_n$ in the expression for $\sigma_n$ is smaller than its exponent in the expression for $U_n$. It is also easy to see, using $\lambda_n/d_n \to 0$ in the case $\beta < 1/2$ and $d_n > Cn^\beta$ when $\beta = 1/2$ or $\beta > 1/2$, that, eventually,

$$\sqrt{n}\sigma_n \ge U_n \sqrt{\log(U_n/\sigma_n)} \asymp U_n \sqrt{\log n}.$$

Then inequality (2.16) gives that

$$E\left( \sup_{t \in C_n} \left| \frac{\Psi(t) \sum_{i=1}^n (K((X_i-t)/h_n) - EK((X-t)/h_n))}{d_n} \right| \right) \le \frac{C\sqrt{n\log n}\,\sigma_n}{d_n}$$

for some $C < \infty$ independent of $n$, as long as $n$ is large enough. For $\beta \le 1/2$, this bound is, up to a multiplicative constant, of the order of

$$\frac{\lambda_n}{d_n} \to 0,$$

and for $\beta > 1/2$, it is of the order of

$$\varepsilon_n^{-(2\beta-1)/2} \sqrt{\log n}\,(n^\beta/d_n)^{1/(2(1-\beta))} h_n^{d\beta/(2(1-\beta))} \to 0,$$

since $d_n \ge Cn^\beta$ for some $C > 0$, and $h_n \to 0$ at least as a negative power of $n$. This completes the proof of the theorem for $\beta < 1$.

For $\beta = 1$, since $d_n \ge Cn$ and $\|\Psi f\|_{\Psi,\infty} \le 1$, we can ignore the centering for all $t$. Then we decompose $B_f$ into $A_n$ defined as above and $B_n := \{t \in B_f : \Psi(t) \le c_n^\beta d_n\}$. The proof of (2.8) and (2.10) with $\lambda_n$ replaced by $d_n$ follows as in the proof of Theorem 2.1, even with some simplification as $B_n$ is now a simpler set. $\square$

We have assumed $d_n \geq Cn^\beta$ and $\beta \leq 1$ in the above theorem. Next we show that these two assumptions are optimal.

REMARK 3.2.   Take $\Psi = f^{-\beta}$. For the sequence (3.1) to be stochastically bounded, it is necessary, by the first part of Theorem 2.1, that the sequence $\{\max_{1\leq i\leq n}(d_n f^\beta(X_i))^{-1}\}$ be stochastically bounded, hence, by regular variation of $d_t$, that

$$\sup n \Pr\left\{\frac{1}{f(X)} \geq d_n^{1/\beta}\right\} < \infty.$$

But if $B_f = \mathbf{R}^d$, then condition (D.c) implies, as in the proof of Corollary 2.2, that

$$n \Pr\left\{\frac{1}{f(X)} \geq d_n^{1/\beta}\right\} \geq c\frac{n}{d_n^{1/\beta}}$$

for all $n$ and some $c > 0$ independent of $n$. Hence, if $\Psi(t)$ is of the order of $f^{-\beta}(t)$, then we must have $d_n \geq Cn^\beta$ in Theorem 3.1.

REMARK 3.3.   Suppose we take $\beta > 1$ in Theorem 3.1, and, again, let us take $\Psi = f^{-\beta}$. Then, we still have that (3.2) is necessary for stochastic boundedness of the sequence (3.1). But then (3.2) implies that

$$\lim_{n\to\infty}\sup_{t\in A_n}\frac{\sum_{i=1}^n K((X_i-t)/h_n)}{d_n f^\beta(t)} = 0 \qquad \text{in pr.}$$

as before. On the other hand, if $B_f = \mathbf{R}^d$, then the set $A_n$ contains $t$'s with $f(t)$ arbitrarily small, and therefore, by (D.b), for some $0 < \delta < \beta - 1$,

$$\sup_{t\in A_n}\frac{nEK((X-t)/h_n)}{d_n f^\beta(t)} \gtrsim \sup_{t\in A_n}\frac{nh_n}{d_n f^{\beta-1-\delta}(t)} = \infty.$$

Hence, the sequence (3.1) is not stochastically bounded, which is a contradiction. So, Theorem 3.1 is not true for $\beta > 1$.

The next theorem describes the almost sure behavior of $\|f_n - Ef_n\|_{\Psi,\infty}$ for large normings.

THEOREM 3.4.   *Assume the usual hypotheses, with condition* (WD.a)$_\beta$ *holding for some* $\beta \in (0,1]$, *and, moreover, that either* $B_f = \mathbf{R}^d$ *or* $K(0) =$ $\|K\|_\infty$. *Let $d_t$ be a strictly increasing regularly varying function satisfying that $\lim_{t\to\infty} d_t/\lambda_t = \infty$ and $d_t \geq Ct^\beta$ for some $C > 0$. Then, either*

$$(3.5) \quad \lim_{n\to\infty} \left\| \frac{\sum_{i=1}^n (K((X_i - t)/h_n) - EK((X-t)/h_n))}{d_n} \right\|_{\Psi,\infty} = 0 \qquad a.s.$$

*or*

$$(3.6) \quad \limsup_{n\to\infty} \left\| \frac{\sum_{i=1}^n (K((X_i - t)/h_n) - EK((X-t)/h_n))}{d_n} \right\|_{\Psi,\infty} = \infty,$$

*according to whether*

$$(3.7) \quad \int_1^\infty \Pr\{\Psi(X) > d_t\}\, dt < \infty \quad or \quad \int_1^\infty \Pr\{\Psi(X) > d_t\}\, dt = \infty.$$

PROOF. Necessity and the part of sufficiency dealing with the sets $A_n$ and $B_n$ follow by a straightforward combination of the proofs of Theorems 2.7 and 3.1. The only difference with previous proofs is in the estimation of the supremum of the processes over the sets

$$C_n = \{t \in B_f : f(t)\Psi(t) \geq \varepsilon_n^{1-\beta} d_n/(nh_n^d)\}.$$

Here, as in the corresponding part of the proof of Theorem 2.6, we use Talagrand's inequality. However, $d_n$ is large and it may fall out of the "Gaussian range" of the inequality. With the notation put forward above, and with the assumptions

$$0 < \sigma_n < U_n/2 \quad \text{and} \quad \sqrt{n}\sigma_n > U_n \sqrt{\log \frac{U_n}{\sigma_n}}$$

shown to hold for all $n$ large enough in the previous proof, Talagrand's inequality in the version from Giné and Guillou [[2001], Proposition 2.2] gives

$$(3.8) \quad \begin{aligned} &\Pr\left\{ \sup_{t\in C_n} \left| \Psi(t) \sum_{i=1}^n \left( K\left(\frac{X_i - t}{h_n}\right) - EK\left(\frac{X-t}{h_n}\right) \right) \right| > \varepsilon d_n \right\} \\ &\leq L \exp\left[ -\frac{1}{L} \frac{\varepsilon d_n}{U_n} \log\left( 1 + \frac{\varepsilon d_n U_n}{Ln\sigma_n^2} \right) \right] := (I), \end{aligned}$$

for some $L$ that depends only on $A$ and $v$ [from (2.15)], and for all $n$ large enough, as long as

$$\frac{\varepsilon d_n}{\sqrt{n}\sigma_n \sqrt{\log U_n/\sigma_n}} > C$$

for a certain constant $C < \infty$. This last condition is eventually satisfied by all $\varepsilon > 0$ since $\log(U_n/\sigma_n) \asymp \log n$ and $d_n/(\sqrt{n \log n}\sigma_n) \to \infty$, as can be easily seen directly from the definitions and properties of these quantities.

Now, by the hypotheses on $h_n$ and $d_n$, and since $\varepsilon_n = 1/\log n$, there exists $\delta > 0$ such that

$$(3.9) \qquad \frac{\varepsilon d_n}{U_n} = \varepsilon \cdot \varepsilon_n^\beta \left(\frac{d_n}{n^\beta}\right)^{1/(1-\beta)} h_n^{-d\beta/(1-\beta)} \geq C\varepsilon \cdot \varepsilon_n^\beta h_n^{-d\beta/(1-\beta)} \geq n^\delta.$$

If $1/2 \leq \beta \leq 1$, then

$$\log\left(1 + \frac{\varepsilon d_n U_n}{Ln\sigma_n^2}\right) \asymp \log\left(1 + \frac{\varepsilon}{L\varepsilon_n^\beta}\right) \gtrsim \log\log n.$$

If $\beta < 1/2$,

$$\log\left(1 + \frac{\varepsilon d_n U_n}{Ln\sigma_n^2}\right) \asymp \log\left(1 + \frac{\varepsilon}{L\varepsilon_n^\beta}\left(\frac{d_n}{nh_n^d}\right)^{(1-2\beta)/(1-\beta)}\right),$$

which is of the order of $\log n$ if the exponent of regular variation of $d_n$ is strictly larger than that of $nh_n^d$, and satisfies

$$\lim_{n\to\infty} n^\delta \log\left(1 + \frac{\varepsilon}{L\varepsilon_n^\beta}\left(\frac{d_n}{nh_n^d}\right)^{(1-2\beta)/(1-\beta)}\right) = \infty$$

for all $\delta > 0$ if the exponents of $d_n$ and $nh_n^d$ coincide. (This can be readily seen using the properties of regular variation and that $\log(1 + \tau) \simeq \tau$ for $\tau$ small.) Combining the last three estimates with the bound (3.9), we get that, for the cases considered,

$$(3.10) \qquad\qquad (I) \leq \exp(-n^\delta)$$

for some $\delta > 0$. Finally, if $\beta < 1/2$ and the exponent of variation of $d_n$ is smaller than the exponent of $nh_n^d$, then

$$\log\left(1 + \frac{\varepsilon d_n U_n}{Ln\sigma_n^2}\right) \simeq \frac{\varepsilon d_n U_n}{Ln\sigma_n^2},$$

and we have, for constants $L$ independent of $n$ (as long as $n$ is large enough) and that vary on each occurrence,

$$
\begin{aligned}
(3.11) \qquad (I) &\leq L\exp\left(-\frac{1}{L}\frac{\varepsilon^2 d_n^2}{n\sigma_n^2}\right) = L\exp\left(-\frac{1}{L}\frac{\varepsilon^2 d_n^2}{nh_n^d}\right) \\
&= L\exp\left(-\frac{1}{L}\varepsilon^2\left(\frac{d_n}{\lambda_n}\right)^2 |\log h_n|\right) \leq L\exp(-M\log n),
\end{aligned}
$$

where $M$ can be made as large as we wish, as long as we take $n$ large enough. (Here we have used $d_n/\lambda_n \to \infty$ and $|\log h_n| \asymp \log n$.) This covers all the cases, and we obtain, combining (3.8), (3.10) and (3.11), that

$$\sum_n \Pr\left\{\sup_{t \in C_n}\left|\Psi(t)\sum_{i=1}^n \left(K\left(\frac{X_i - t}{h_n}\right) - EK\left(\frac{X - t}{h_n}\right)\right)\right| > \varepsilon d_n\right\} < \infty$$

for all $\varepsilon > 0$, proving that

$$\lim_{n \to \infty} \sup_{t \in C_n} \left| \Psi(t) \sum_{i=1}^{n} \left( K\left( \frac{X_i - t}{h_n} \right) - EK\left( \frac{X - t}{h_n} \right) \right) \right| = 0 \quad \text{a.s.}$$

This completes the proof of the theorem. $\square$

The results in this section obviously apply to the densities in Examples 2.8–2.11. For instance, let $f$ be the symmetric exponential density on $\mathbf{R}$ considered in Examples 2.8 and 2.10, and let $h_n = n^{-\alpha}$, $0 < \alpha < 1$. Then, Theorem 3.1 shows that

$$n^{1-\alpha-\beta} \left\| \frac{f_n - Ef_n}{f^\beta} \right\|_\infty \xrightarrow{d} \|K\|_\infty Z^\beta,$$

where $Z$ is the random variable defined in Example 2.10, if and only if

$$\frac{1-\alpha}{2} < \beta < 1;$$

and Theorem 3.4 shows that, for $c(t)$ strictly increasing and regularly varying,

$$\frac{n^{1-\alpha}}{c^\beta(n)} \left\| \frac{f_n - Ef_n}{f^\beta} \right\|_\infty \to 0 \qquad \text{a.s.}$$

if and only if

$$\int^\infty \frac{dt}{c(t)} < \infty.$$

A similar statement holds true for normal densities.

REMARK 3.5. Suppose that $\mathcal{K}$ is a uniformly bounded class of kernels supported by a fixed bounded set and such that the class

$$\mathcal{F} := \left\{ K\left( \frac{\cdot - t}{h} \right) : t \in \mathbf{R}^d, \ h > 0, \ K \in \mathcal{K} \right\}$$

is measurable and has covering numbers

$$N(\mathcal{F}, L_2(P), \|K\|_{L_2(P)}\varepsilon) \le \left( \frac{A}{\varepsilon} \right)^v, \qquad 0 < \varepsilon < 1,$$

for some $A$ and $v$ finite and positive and for all probability measures $P$. [In particular, $\mathcal{K}$ may be a subset of the linear span of a finite set of functions $k$ as defined in condition (K)]. Suppose we wish to consider

$$\sup_{K \in \mathcal{K}} \left\| \frac{f_n - Ef_n}{c_n} \right\|_{\Psi, \infty},$$

where $c_n$ is $d_n$ or $\lambda_n$, as defined above. Then uniform boundedness and uniformity of the support allow us to deal with the sup over $A_n$ and $B_n$, and the entropy bound, with the sup over $C_n$, just as in the previous theorems. The sup over the central part $D_a$ is handled in Mason (2004). So, it is straightforward to prove a uniform in $K \in \mathcal{K}$ version of our results. It is also possible to prove a functional law of the logarithm in our setting by following Mason (2004).

**Acknowledgments.** We thank David M. Mason for several comments on previous drafts and for useful conversations on the subject of this article. Lyudmila Sakhanenko read the manuscript and pointed out several typos. This research has been carried out, aside from our institutions of origin, at the Department of Mathematics of the University of Connecticut (V. Koltchinskii and J. Zinn), at the Department of Statistics of the University of Washington, Seattle (V. Koltchinskii), and at the Institut Henri Poincaré, Paris (E. Giné and V. Koltchinskii). It is a pleasure to thank these institutions for their the support and hospitality.

## REFERENCES

Deheuvels, P. (2000). Uniform limit laws for kernel density estimators on possibly unbounded intervals. In *Recent Advances in Reliability Theory: Methodology, Practice and Inference* (N. Limnios and M. Nikulin, eds.) 477–492. Birkhäuser, Boston.

de la Peña, V. and Giné, E. (1999). *Decoupling, from Dependence to Independence.* Springer, New York. MR1783500

Dudley, R. M. (1999). *Uniform Central Limit Theorems.* Cambridge Univ. Press. MR1720712

Einmahl, J. H. J. and Mason, D. (1985a). Bounds for weighted multivariate empirical distribution functions. *Z. Wahrsch. Verw. Gebiete* **70** 563–571. MR807337

Einmahl, J. H. J. and Mason, D. (1985b). Laws of the iterated logarithm in the tails for weighted uniform empirical processes *Ann. Probab.* **16** 126–141. MR920259

Einmahl, J. H. J. and Mason, D. (1988). Strong limit theorems for weighted quantile processes. *Ann. Probab.* **16** 1623–1643. MR958207

Einmahl, U. and Mason, D. (2000). An empirical process approach to the uniform consistency of kernel-type function estimators. *J. Theoret. Probab.* **13** 1–37. MR1744994

Feller, W. (1971). *An Introduction to Probability Theory and Its Applications* **2**. Wiley, New York. MR270403

Giné, E. and Guillou, A. (2001). On consistency of kernel density estimators for randomly censored data: Rates holding uniformly over adaptive intervals. *Ann. Inst. H. Poincaré Probab. Statist.* **37** 503–522. MR1876841

Giné, E. and Guillou, A. (2002). Rates of strong consistency for multivariate kernel density estimators. *Ann. Inst. H. Poincaré Probab. Statist.* **38** 907–922. MR1955344

Giné, E. and Zinn, J. (1984). Some limit theorems for empirical processes. *Ann. Probab.* **12** 929–989.

Koltchinskii, V. and Sakhanenko, L. (2000). Testing for ellipsoidal symmetry of a multivariate distribution. In *High Dimensional Probability II* (E. Giné, D. M. Mason and J. A. Wellner, eds.) 493–510. Birkhäuser, Boston. MR757767

MASON, D. (2004). A uniform functional law of logarithm for the local empirical process. *Ann. Probab.* **32** 1391–1418. MR2060302

MONTGOMERY-SMITH, S. J. (1993). Comparison of sums of independent identically distributed random vectors. *Probab. Math. Statist.* **14** 281–285. MR1321767

NOLAN, D. and POLLARD, D. (1987). *U*-processes: Rates of convergence. *Ann. Statist.* **15** 780–799. MR888439

PARZEN, E. (1962). On the estimation of a probability density function and the mode. *Ann. Math. Statist.* **33** 1065–1076. MR143282

POLLARD, D. (1984). *Convergence of Stochastic Processes.* Springer, New York. MR762984

ROSENBLATT, M. (1956). Remarks on some nonparametric estimates of a density function. *Ann. Math. Statist.* **27** 832–835. MR79873

SILVERMAN, B. W. (1978). Weak and strong uniform consistency of the kernel estimate of a density and its derivatives. *Ann. Statist.* **6** 177–189. MR471166

STOUT, W. (1974). *Almost Sure Convergence.* Academic Press, New York. MR455094

STUTE, W. (1982). A law of the logarithm for kernel density estimators. *Ann. Probab.* **10** 414–422. MR647513

STUTE, W. (1984). The oscillation behavior of empirical processes: The multivariate case. *Ann. Probab.* **12** 361–379. MR735843

TALAGRAND, M. (1994). Sharper bounds for Gaussian and empirical processes. *Ann. Probab.* **22** 28–76. MR1258865

TALAGRAND, M. (1996). New concentration inequalities in product spaces. *Invent. Math.* **126** 505–563. MR1419006

E. GINÉ
DEPARTMENT OF MATHEMATICS
DEPARTMENT OF STATISTICS
UNIVERSITY OF CONNECTICUT
STORRS, CONNECTICUT 06269
USA
E-MAIL: gine@uconnvm.uconn.edu

V. KOLTCHINSKII
DEPARTMENT OF MATHEMATICS
AND STATISTICS
UNIVERSITY OF NEW MEXICO
ALBUQUERQUE, NEW MEXICO 87131
USA
E-MAIL: vlad@math.unm.edu

J. ZINN
DEPARTMENT OF MATHEMATICS
TEXAS A&M UNIVERSITY
COLLEGE STATION, TEXAS 77843
USA
E-MAIL: jzinn@math.tamu.edu